\documentstyle[12pt,amscd,amsfonts]{amsart}
\newtheorem{thm}[equation]{Theorem}
\numberwithin{equation}{section}
\newtheorem{cor}[equation]{Corollary}

\newtheorem{diag}[equation]{Diagram}

\newtheorem{lem}[equation]{Lemma}
\newtheorem{conj}[equation]{Conjecture}
\newtheorem{defin}[equation]{Definition}
\newtheorem{prop}[equation]{Proposition}

\newsymbol\square 1003
\begin{document}
\raggedbottom
\voffset=-.7truein
\hoffset=0truein
\vsize=8truein
\hsize=6truein
\textheight=8truein
\textwidth=6truein
\baselineskip=18truept
\def\mapright#1{\smash{\mathop{\longrightarrow}\limits^{#1}}}
\def\mapleft#1{\smash{\mathop{\longleftarrow}\limits^{#1}}}
\def\mapup#1{\Big\uparrow\rlap{$\vcenter {\hbox {$#1$}}$}}
\def\mapdown#1{\Big\downarrow\rlap{$\vcenter {\hbox {$\ssize{#1}$}}$}}
\def\mapne#1{\nearrow\rlap{$\vcenter {\hbox {$#1$}}$}}
\def\mapse#1{\searrow\rlap{$\vcenter {\hbox {$\ssize{#1}$}}$}}
\def\mapr#1{\smash{\mathop{\rightarrow}\limits^{#1}}}
\def\lb{[}
\def\ss{\smallskip}
\def\sm{\wedge}
\def\la{\langle}
\def\ra{\rangle}
\def\on{\operatorname}
\def\kbar{{\overline k}}
\def\qed{\quad\rule{8pt}{8pt}\bigskip}
\def\ssize{\scriptstyle}
\def\a{\alpha}
\def\bz{{\bold Z}}
\def\fr{\on{Fr}}
\def\im{\on{im}}
\def\ext{\on{Ext}}
\def\cm{{\cal M}}
\def\sq{\on{Sq}}
\def\eps{\epsilon}
\def\ar#1{\stackrel {#1}{\rightarrow}}
\def\br{\bold R}
\def\bc{\bold C}
\def\si{\sigma}
\def\Ebar{{\overline E}}
\def\Sum{\sum}
\def\tfrac{\textstyle\frac}
\def\tb{\textstyle\binom}
\def\Si{\Sigma}
\def\w{\wedge}
\def\equ{\begin{equation}}
\def\b{\beta}
\def\G{\Gamma}
\def\g{\gamma}
\def\endeq{\end{equation}}
\def\sn{S^{2n+1}}
\def\zp{\bold Z_p}
\def\P{{\cal P}}
\def\cG{{\cal G}}
\def\Hom{\on{Hom}}
\def\ker{\on{ker}}
\def\coker{\on{coker}}
\def\da{\downarrow}
\def\io{\iota}
\def\Om{\Omega}
\def\u{{\cal V}}
\def\e{{\cal E}}
\def\exp{\on{exp}}
\def\xbar{{\overline x}}
\def\ebar{{\overline e}}
\def\et{{\widetilde E}}
\def\ni{\noindent}
\def\coef{\on{coef}}
\def\den{\on{den}}
\def\lcm{\on{l.c.m.}}
\def\vi{v_1^{-1}}
\def\ot{\otimes}
\def\psibar{{\overline\psi}}
\def\phibar{{\overline\phi}}
\def\mhat{{\hat m}}
\def\exc{\on{exc}}
\def\ms{\medskip}
\def\ehat{{\hat e}}
\def\dirlim{\on{dirlim}}
\def\psit{\widetilde\psi}
\def\bq{{\bold Q}}
\def\vp{v_1^{-1}\pi}
\title[3-primary periodic homotopy groups of $E_7$]
{3-primary $v_1$-periodic homotopy groups of $E_7$}
\author[D. Davis]{Donald M. Davis}
\address{Lehigh University\\Bethlehem, PA 18015}
\email{dmd1@@lehigh.edu}

\date{}

\maketitle
\section{Introduction}\label{intro}
In this paper we compute the 3-primary $v_1$-periodic homotopy groups of the
exceptional Lie group $E_7$.

The $p$-primary $v_1$-periodic homotopy groups of a space $X$, denoted
$\vp_*(X;p)$ or just $\vp_*(X)$,
were defined in \cite{DM}. They are a localization of the actual
homotopy groups, telling roughly the portion which is detected by $K$-theory
and its operations. If $X$ is a compact Lie group, each $\vp_i(X;p)$ is a
direct summand of some actual homotopy group of $X$, and so summands of
$v_1$-periodic homotopy groups of $X$ give lower bounds for the $p$-exponent of
$X$.

After the author computed $\vp_*(SU(n);p)$ for odd $p$ in 1989, Mimura proposed
the goal of calculating $\vp_*(X;p)$ for all compact simple Lie groups $X$.
This has now been achieved in the following cases $(X,p)$:
\begin{itemize}
\item $X$ a classical group and $p$ odd (\cite{D});
\item $X$ an exceptional Lie group with $H_*(X;\bold Z)$ $p$-torsion-free
(\cite{BDMi});
\item ($SU(n)$ or $Sp(n)$, $2$) (\cite{BDSU},\cite{BDM});
\item $(G_2,2)$ (\cite{DMG}), ($F_4$ or $E_6$, 3) (\cite{F4}),
and $(E_7,3)$ (the current paper).
\end{itemize}

The only cases remaining then are $(E_8,2\text{ or }3\text{ or }5)$ and $(SO(n)
\text{ or } F_4\text{ or } E_6\text{ or }E_7,2)$. Several of these appear
tractable.

Now we state our main theorem. We usually abbreviate $\vp_*(X;3)$ as
$v_*(X)$, and denote by $\nu(n)$ the exponent of 3 in the integer $n$.
\begin{thm} \label{main}
If $j$ is even, then $v_{2j}(E_7)=v_{2j-1}(E_7)=0$. If $j$ is odd, then
$$v_{2j}(E_7)\approx v_{2j-1}(E_7)\approx\begin{cases}
\bz/3\oplus\bz/3^{\min(10,\nu(j-9-2\cdot3^5)+4)}&\text{if $j\equiv0$ mod $3$}\\
\bz/3\oplus\bz/3^{\min(8,\nu(j-43)+5)}&\text{if $j\equiv1,7$ mod $9$}\\
\bz/3\oplus\bz/3^{\min(14,\nu(j-13-4\cdot3^8)+5)}&\text{if $j\equiv4$ mod $9$}\\
\bz/9\oplus\bz/3^{\min(19,\nu(j-17-2\delta\cdot3^{13})+4)}&\text{if $j\equiv5,8$
mod $9$,}\end{cases}$$
where $\delta$ equals one of the numbers $2$, $5$, or $8$.
If $j$ is odd and $j\equiv2$ mod $9$,
then $v_{2j-1}(E_7)\approx\bz/9\oplus\bz/3^{\min(13,\nu(j-11)+4)}$, while
$$v_{2j}(E_7)\approx\begin{cases}\bz/3^3\oplus\bz/3^{\nu(j-11)+3}&
\text{if $\nu(j-11)<10$}\\
\bz/3^3\oplus\bz/3^{12}\text{ or }\bz/3^4\oplus\bz/3^{11}&\text{if
$\nu(j-11)\ge10$.}\end{cases}$$
\end{thm}

An immediate corollary of this work is a lower bound for the 3-exponent of
$E_7$. Recall that the $p$-exponent of a space $X$, denoted $\exp_p(X)$, is the
largest $e$ such that $\pi_*(X)$ has an element of order $p^e$. We obtain
\begin{cor} The $3$-exponent of $E_7$ satisfies $\exp_3(E_7)\ge19$.
\end{cor}
\begin{pf} If $X$ is a compact Lie group, then $v_i(X)\approx\dirlim_{k,m}
\pi_{i+4k\cdot3^m}(X)$. Hence an element of order $3^{19}$ in $v_{2j}(E_7)$
when $j\equiv17+2\delta\cdot3^{13}$ mod $2\cdot3^{15}$
corresponds to an element of order $3^{19}$ in some $\pi_n(E_7)$.
\end{pf}
By comparison, the result that we have obtained at other primes is (\cite{BDMi})
$$\exp_p(E_7)\begin{cases}=17&\text{if $p>17$}\\
=18&\text{if $p=17$}\\
\ge17&\text{if $p=7$, 11, or 13}\\
\ge18&\text{if $p=5$.}\end{cases}$$
This should be contrasted with the situation for spheres, where we have
$\exp_p(\sn)=n$ for all odd primes $p$ and all positive integers $n$ by
\cite{CMN}.

Note that in Theorem \ref{main}, we determine the precise abelian group
structure of all groups (with isolated exceptions),
whereas in some earlier papers, such as
\cite{D}, \cite{BDMi}, and \cite{F4}, we had been unable to determine the
group structure of most groups $v_{2j-1}(X)$. Because of the insights of
\cite{DY}, we are able not only to resolve the extension questions
(group structure) in almost all cases occurring here,
but also in those of
\cite{BDMi} and \cite{F4}. These new results about group structure
are presented in Section \ref{extensions}.

Most of the work is calculation of the $v_1$-periodic
unstable Novikov spectral sequence
(UNSS) of the space $Y_7:=\Om E_7/Sp(2)$.
The main input is the detailed structure of $H_*(\Om E_7;\bz/3)$ given in
\cite{HH} and restated here in Proposition \ref{HaHa}. The advantage of
$Y_7$ over $\Om E_7$ is that $BP_*(Y_7)$ is a free commutative algebra,
which makes its UNSS easier to calculate. Perhaps the most
novel feature of the calculations here is the use of coassociativity
to give detailed formulas for $BP_*$-coaction. The terms which arise in this
way play crucial roles in the calculations. The calculations of $v_*(F_4)$
in \cite{F4} are essential in the transition from $v_*(Y_7)$ to $v_*(\Om E_7)$.

Another delicate point is convergence of the $v_1$-periodic UNSS for $Y_7$.
In Section \ref{E7F4}, we use deep recent work of Bousfield and
Bendersky-Thompson to prove that the $v_1$-periodic UNSS converges
to $v_*(-)$ for
$E_7/F_4$, which we will show implies similar convergence for $Y_7$.

The author would like to thank Mamoru Mimura, Pete Bousfield, and
especially Martin Bendersky for useful suggestions.

\section{Background in $v_1$-periodic homotopy and the UNSS}\label{background}

In this section, we review known results, and establish one new useful result
about computing the UNSS.
Although some of these results are
also true when $p=2$, it will simplify exposition to assume that $p$ is an odd
prime.

The $v_1$-periodic homotopy groups of any topological space $X$ are defined
 by
\begin{equation}\label{vdef}\vp_i(X)=\lim_{\to}[M^{i+1+kqp^e}(p^e),X],
\end{equation}
where $q=2p-2$, a notation that will be used consistently throughout this
paper, and $M^t(n)$ denotes the Moore space $S^{t-1}\cup_ne^t$.
Here the direct limit is taken over increasing values of $e$ and $k$
using Adams maps $M^{t+qp^e}(p^e)\to M^t(p^e)$ and canonical maps
$M^t(p^{e+1})\to M^t(p^e)$.
This definition was given in \cite{DM}, where their relationship
with actual homotopy groups of many spaces was established.

A space $X$ is said to {\em have an $H$-space exponent} at the prime $p$
if, for some $e$ and $L$, $p^e:\Omega^L X\to\Omega^L X$ is null homotopic.
It was shown in \cite[1.9]{DM} that if $X$ has an $H$-space exponent, then
$$\vp_i(X)\approx\lim_{\to}\pi_{i+kqp^e}(X),$$
and hence $\vp_i(X)$ is a direct summand of some group $\pi_{i+kqp^e}(X)$.
To make this final deduction, we need to know that the limit group is a finitely
generated abelian group, but this will be the case.

Next we discuss the unstable cobar complex, which can be used to compute the
UNSS for many spaces.  We will modify and
generalize previous treatments of this topic.
Let $BP$ be the Brown-Peterson spectrum corresponding to the prime $p$. Then
$$BP_*=\pi_*(BP)\approx Z_{(p)}[v_1, v_2,\ldots ],$$
where $v_i$ are  the Hazewinkel generators of $BP_*$.
Let $\Gamma=BP_*(BP)\approx BP_*[h_1,h_2,\ldots ]$, where $h_i$ are
conjugates of Quillen's
generators $t_i$. We have $|v_i|=|h_i|=2(p^i-1)$.
Let $\eta=\eta_R: BP_*\to BP_*(BP)$
be the right unit. We write $h_iv_j$ interchangeably with
$\eta(v_j)h_i$; this is the right action of $BP_*$ on $\Gamma$.

     Let $M$ be a $\Gamma$-comodule with coaction map
$\psi_M: M\mapr{} \Gamma\otimes M$. Tensor products are always over $BP_*$.
The stable cobar complex $SC^*(M)$ is defined by
$$SC^s(M)=\Gamma\otimes\Gamma\otimes\cdots \otimes
\Gamma\otimes M,$$
with $s$ copies of $\Gamma$,
and differential $d:SC^s(M)\to SC^{s+1}(M)$ given by
\begin{eqnarray}
d(\g_1\ot\cdots\ot\g_s\ot m)&=&1\ot\g_1\ot\cdots\g_s\ot m\\
&&+
\sum_{j=1}^s (-1)^j\g_1\ot\cdots\ot\psi(\g_j)\ot\cdots\ot\g_s\ot m\nonumber\\
&&+(-1)^{s+1}\g_1\ot\cdots\ot\g_s\ot\psi_M(m).\label{d}
 \end{eqnarray}

Our unstable cobar complex $VC^*(M)$ is
 a subcomplex of $SC^*(M)$, consisting of terms satisfying an unstable
condition, introduced in the following definition.

\begin{defin}{\rm \cite[3.3]{Der}} \label{unstable}
If $M$ is a nonnegatively graded free left $A$-module,
 then $V(M)$ is defined to be the $BP_*$-span of
$$\{h^I\otimes m\ :\  2(i_1+i_2+\cdots )\le |m|\}\subset \Gamma\otimes
 M,$$
where $I=(i_1, i_2, \ldots )$ and $h^I=h_1^{i_1}h_2^{i_2}\cdots $.\end{defin}
\noindent This unstable condition will pervade our computations. Note that for
odd-dimensional classes, this agrees with the module $U(M)$ which has been used
most frequently in earlier work of the author and Bendersky. However, it
also agrees with the $V(M)$ construction employed in \cite{BDM}
on even-dimensional classes. The novelty
here is that it will be applied to a module having classes of both parities.

Define $VC^0(M)=M$, and
$VC^s(M)=V(VC^{s-1}(M))$. If $M$ is a $\Gamma$-comodule,
then the differential $d$ of the stable cobar complex of $M$
induces a differential on the subcomplex $VC^*(M)$.
 We will usually replace it by the
chain-equivalent reduced complex obtained by replacing $V(M)$ by
$\ker(V(M)\mapright{\eps} M)$.(\cite[2.16]{RIMS})
 This has the effect of only looking at terms
which have positive grading in each position.
The homology groups of this unstable cobar complex are denoted
by $\ext_{\u}^{s,t}(M)$. As observed in \cite{Der}, these are the usual Ext
groups in the abelian category $\u$ of $\Gamma$-comodules satisfying the
unstable condition in Definition \ref{unstable}.
Note there is a shift isomorphism
\begin{equation}\label{shift}\ext_{\u}^{s,t-1}(BP_*S^{2n})
\approx\ext_{\u}^{s,t}(BP_*S^{2n+1}),\end{equation}
induced by a shift isomorphism of the unstable cobar complexes.

The following generalization of \cite[\S7]{BCM} will be very useful to us.
Its proof follows some suggestions of Martin Bendersky.
\begin{thm}\label{Vthm}If $X$ is a simply-connected $CW$-space,
 there is a spectral sequence $\{E^{s,t}_r(X), d_r\}$
which converges to the homotopy groups of $X$ localized at $p$. If $X$ is an
$H$-space, and
 $BP_*(X)$ is a
free commutative algebra, then
              $$E_2^{s,t}(X)=\ext^{s,t}_{\u}(Q(BP_* X)),$$
where $Q(BP_*X)$ denotes the indecomposable quotient
of $BP_*X$.\end{thm}
This is the UNSS for the space $X$. We will write $VC^*(X)$ for the complex
$VC^*(Q(BP_*X))$, whose homology is $E_2(X)$.
We denote by $\fr$  the free commutative algebra functor. If $N$ is
a free $BP_*$-module
with basis $B=B_{\text{ev}}\cup B_{\text{od}}$,
then $\fr(N)$ is the tensor product of a polynomial algebra over
$BP_*$ on $B_{\text{ev}}$  with an exterior algebra
on $B_{\text{od}}$.

\begin{pf} The spectral sequence was described in \cite{BCM}.
The determination of $E_2$ when $BP_*(X)$ is a free commutative algebra
is quite similar to that of \cite[6.1]{BCR} and to the
argument on \cite[p.346]{BDM}. Let $M=Q(BP_*X)$, a free $BP_*$-module.

Let $\cal G$ denote the category of unstable $\G$-coalgebras, and $G(-)$
the associated functor considered in \cite{BCM}. If $N$ is a free
$BP_*$-module, then $G(N)$ is defined to be $BP_*(BP(N))$, where $BP(N)$
is the 0th space of the $\Om$-spectrum representing the homology theory
$BP_*(-)\ot N$.
If $N$ has basis $B$, then
\begin{equation}\label{Gstr}
G(N)\approx BP_*(\prod_{b\in B}\bold{BP}_{|b|})\approx
\fr(\langle h^Ib:b\in B,\  2|I|\le|b|\rangle).\end{equation}
Here $h^I$ is as in Definition \ref{unstable} with $|I|=\sum i_j$,
while $|b|$ denotes the degree of the basis element
$b$. Also $\bold{BP}_n$ denotes the $n$th space in the $\Om$-spectrum for $BP$.
The first isomorphism in (\ref{Gstr}) is immediate from the definition of $G$
given in \cite[6.3,6.7]{BCM}. The second isomorphism follows from
\cite[p.51]{Wil}, which says that $H_*(\bold{BP}_n)$ is a polynomial algebra if
$n$ is even, and an exterior algebra if $n$ is odd, \cite[4.9]{HR}, which says
that the same thing is then
true of $BP_*(\bold{BP}_n)$, and \cite[p.1040]{BPH},
which interprets conveniently
the description of the indecomposables first given in \cite{HR}.
Note that there is an isomorphism of $BP_*$-modules
\begin{equation}\label{QGV}Q(G(N))\approx V(N).\end{equation}

We claim that
\begin{equation}
 BP_*X\mapright{\xi} G(M)@>\to>> G(V(M))@>\to>\to>G(V^2(M))\cdots
\label{Gseq}\end{equation}
is an augmented cosimplicial resolution in $\cal G$.
Here the augmentation $\xi$ is the composite
$$BP_*X\mapright{\eta_X}G(BP_*X)@>G(\rho)>> G(QBP_*X),$$
where the second morphism applies $G$ to the quotient morphism $\rho$.
The cofaces  are of two types:
\begin{itemize}
\item
$G(V^qM)\mapright{\eta_{G(V^qM)}}G(G(V^qM))\mapright{G(\rho)}G(V^{q+1}M)$,
where $\rho:G(-)\to QG(-)=V(-)$ is the quotient morphism.
\item $G(V^i(\psi_{V^{q-i}M}))$, $0\le i\le q$, where $\psi_N:N\to V(N)$
stabilizes to the $\G$-coaction.
\end{itemize}
The degeneracies $G(V^qM)\to G(V^{q-1}M)$ just do the counit $\eps$ on one of
the $V$-factors.
It is clear that all of these morphisms are in $\cal G$, and the
cosimplicial identities are satisfied as usual. The argument of
\cite[3.13]{BDM} implies that the first type of coface map and the augmentation
$\xi$ are algebra morphisms. The second type of coface map is an algebra
morphism since it is $BP_*(f)$ for an infinite loop map $f$, namely the map
$BP(N)\mapright{BP(g)}BP(N')$ induced by a $BP_*$-morphism $N\mapright{g}N'$.

The exactness of the resulting augmented
cochain complex
\begin{equation}\label{acy}
0\to BP_*X\mapright{\xi} G(M)\to G(V(M))\to G(V^2(M))\to\cdots\end{equation}
(obtained using the alternating sum of cofaces as boundaries)
follows as in \cite[p.387]{BCR}, but we provide details for
completeness.  (In comparing with \cite{BCR}, it is useful to note
that $V(N)\approx\sigma^{-1}U(\sigma N)$.)
Since the coface operators
are algebra homomorphisms, their alternating sum preserves the filtration
of this augmented complex by
powers of the augmentation ideal.
Let $E_0$ denote the quotients of the filtration.  Then,
using (\ref{QGV}), we have
$$E_0(G(V^q(M)))\approx \fr(Q(G(V^q(M))))\approx \fr(V^{q+1}(M)),$$
and $E_0(BP_*X)\approx \fr(M)$.
Thus $E_0(\ref{acy})$
is the free commutative algebra on the complex
$$0\to M\to V(M)\to V^2(M)\to\cdots$$
with morphisms the alternating sum of $\psi$ on each $V$ and $\psi_M$,
 which is exact by \cite[7.8]{Bous}. Since the free commutative algebra functor
applied to an exact sequence yields an exact sequence, we deduce that
(\ref{acy}) is exact, and hence yields a resolution in $\cal G$ of $BP_*X$.

Hence $\ext_{\cal G}(BP_*,BP_*X)$ is equal to the cohomology
of the complex obtained by applying $\Hom_{\cal G}(BP_*,-)$ to the portion of
(\ref{acy}) after $\xi$. Since
$\Hom_{\cal G}(BP_*,G(N))\approx N$, we obtain that $\ext_{\cal G}(BP_*,BP_*X)$
is the homology of the complex
\begin{equation}\label{UCC}M\to V(M)\to V^2(M)\to\cdots,\end{equation}
with differentials as in (\ref{d}).
The claim of the theorem follows now from \cite[6.17]{BCM},
which states that $E_2(X)\approx \ext_{\cal G}(BP_*,BP_*X)$, and the
observation that (\ref{UCC}) is just our unstable cobar complex, whose homology
is $\ext_\u(M)$.
\end{pf}

The following definition will be extremely important.
\begin{defin}\label{exdef} The excess {\em{exc}}$(\g)$ of an element $\g$ of
$\bar\Gamma^s$ is defined to be the smallest $n$ such that $\g\io_{2n+1}$
is an element of $VC^s(\sn)$.\end{defin}

This means that if $\g=\g_1\ot\cdots\ot\g_s$, then for $1\le i\le s$,
$$\g_i\ot(\g_{i+1}\cdots\g_s\io_{2n+1})$$
must satisfy \ref{unstable}. The following result, which was proved as
\cite[4.2]{exp}, gives a formula for the excess of certain monomials when
$s=2$.
\begin{lem} If $a\le b$ and $a\le d$, then
$$\exc(p^ah^b\ot v^ch^dv^e)=\max\biggl(b-(p-1)(c+d),d\biggr)-\min\biggl(a,
|b-(p-1)c-pd|\biggr)-(p-1)e.$$\label{excfor}
\end{lem}

In \cite{Bloc}, the $v_1$-periodic UNSS was defined and shown to satisfy the
following very nice property.
\begin{thm} \label{perUNSS} If $p$ is odd and $X$ is spherically resolved,
the $v_1$-periodic UNSS of $X$ satisfies
\begin{itemize}
\item $\vi E_\infty^{s,t}(X)=\vi E_2^{s,t}(X)$, and is
$0$ unless $s=1$ or $2$ and $t$ is odd.
\item $\vi E_\infty^{s,t}(X)\approx\vp_{t-s}(X)$ if $s=1$ or $2$ and $t$ is odd.
\item $\vi E_2^{s,t}(X)=\dirlim E_2^{s,t+kqp^e}(X)$, where $e$ is chosen
sufficiently large, and the direct limit is taken over increasing values of $k$
under multiplication by $v_1^{p^e}$.
\end{itemize}
\end{thm}
Here we say that $X$ is spherically resolved if it can be built from a finite
number of odd-dimensional spheres by fibrations. In Section \ref{E7F4}, we will
show that Theorem \ref{perUNSS} holds in a certain case in which we cannot
prove that $X$ is spherically resolved.

We will use the unstable cobar complex for the unlocalized UNSS, but, as we are
dealing exclusively with $v_1$-periodic classes, we can, in effect, act as if
it satisfies the first two properties of Theorem \ref{perUNSS}.

We will make frequent use of the following result for the spheres, which was
proved in \cite{exp}, following \cite{BCM} and \cite{uns}. We introduce here
terminology $x\equiv y$ mod $S^{2n-1}$ to mean that $x-y$ desuspends to (or is
defined on) $S^{2n-1}$. For elements of $E_2^s(\sn)$, we frequently abbreviate
$x\iota_{2n+1}$ as $x$.

\begin{thm} \label{sph}
\begin{enumerate}
\item The only nonzero groups $\vi E_2^{s,t}(\sn)$ are
$$\vi E_2^{s,2n+1+qm}(\sn)\approx\bz/p^e$$
with $s=1$ or $2$ and $e=\min(n,\nu(m)+1)$.
\item
The generator of
$E_2^{1,2n+1+qm}(\sn)$ is $\a_{m/e}:=d(v_1^m)/p^e$ and satisfies
\begin{equation}\a_{m/e}\equiv -v_1^{m-e}h_1^e \mod S^{2e-
1},\label{e}\end{equation}
and, if $m=sp^{e-1}$ with $s\not\equiv0$ mod $p$, and $e> n$, then
\begin{equation}\label{s}\a_{m/e}\equiv-sv_1^{m-1}h_1\mod p.\end{equation}

\item
If $n\le\nu(m)+1$
and $1\le j\le n$, then $d(p^{m-\nu(m)-1-j}h_1^m)\io_{2n+1}$ has order $p^j$ in
$E_2^{2,2n+1+qm}(\sn)$. It equals
$v_1^{m-j-1}h_1\otimes h_1^j$ mod $S^{2j-1}$.
\item
If $\nu(m)+1\le n$ and $1\le j\le\nu(m)+1$, then $d(p^{m-n-j}h_1^m)\io_{2n+1}$
has order $p^j$ in $E_2^{2,2n+1+qm}(\sn)$. It equals
$v_1^{m-n-j+\nu(m)}h_1\otimes h_1^{n+j-\nu(m)-1}$ mod $S^{2n+2j-2\nu(m)-3}$.
\item The homomorphism $\Sigma^2:E_2^{2,2n-1+qm}(S^{2n-1})\to
E_2^{2,2n+1+qm}(\sn)$ is injective if $n\le\nu(m)+1$ and is multiplication by
$p$ otherwise.
\end{enumerate}
\end{thm}

Other more technical
results proved in earlier works are as follows. Here we
begin the practice, which will be continued throughout the paper, of often
abbreviating $h_1$ as $h$, and $v_1$ as $v$. Also, we introduce the term
\lq\lq leading term'' to refer to a monomial of largest excess in an element
$z$ of $VC(X)$; all other monomials comprising $z$ desuspend farther than does
the leading term.
\begin{prop}\label{sph2}
\begin{enumerate}
\item $($\cite[2.9]{DY}$)$ If a cycle of $VC^{2}(\sn)$ has order $p^f$
in $E_2^{2,t}(\sn)$ and
leading term $h\ot h^j\io_{2n+1}$, then $j+\nu(|E_2^{2,t}(\sn)|)=f+n$.
\item $($\cite[4.6]{exp}$)$ Let $\nu=\nu(\sigma)$, and let
$$z=\eps v^{\sigma-e-1}h^e\ot h+L\in VC^{2,2n+1+q\sigma}(\sn)$$
be a cycle with $\eps\in\bz_{(p)}$ and $\exc(L)<e-p+1\le n-\nu$. Then
$$z=d(u\eps v^{\sigma-(e+\nu-p+2)}h^{e+\nu-p+2}+L'),$$
where $u$ is a unit in $\bz_{(p)}$, and $\exc(L')<e+\nu-p+2$. The same
conclusion holds for $z=\eps v^{\sigma-e+p-2}h\ot h^{e-p+1}+L$.
\end{enumerate}
\end{prop}

We will need the following precise description of $\a_2$.
\begin{lem}\label{a2} The element $\a_2$ which generates $E_2^{1,2q+2n+1}(\sn)$
is given by
$$\a_2=-d(v_1^2)/p=\tfrac1p(v^2-(v-ph)^2)=2vh-ph^2=hv+vh.$$
\end{lem}

We will make repeated use of the following result, especially part (1).
\begin{lem}\label{std} Let $p=3$. Then

\begin{enumerate}
\item $\eta(v_1)=v_1-3h_1$
\item $\eta(v_2)=v_2+4v^3h-18v^2h^2+35vh^3-24h^4-3h_2$
\item $\psi(h_1)=h_1\ot 1+1\ot h_1$
\item $\psi(h_2)=h_2\ot1+1\ot h_2+4h^3\ot h+6h^2\ot h^2+3h\ot h^3-vh\ot
h^2-vh^2\ot h$
\end{enumerate}
\end{lem}
\begin{pf} Parts (1) and (3) are standard, appearing in all referenced papers
of the author and/or Bendersky. Part (2) is taken from Giambalvo's tables
(\cite{Giam}). Part (4) is derived from \cite[2.6i]{BDMi}, using part (1) of
this lemma several times to replace a $v$ on the right (which is interpreted as
$\eta(v)$) by $v-3h$. Note, however, that the sum in \cite[2.6i]{BDMi} should
be preceded by a minus sign.
\end{pf}

The following result, proved in \cite[2.11,2.12,2.13]{DY}, will be central to
many of our calculations.
\begin{lem} \label{yang} \begin{enumerate}
\item If $n\ge1$, then in $E_2^1(\sn)$, $h_1^pv_1\equiv v_1^ph_1$ mod $S^1$.
\item $h_1^{n+p-1}\ot h_1\equiv -v_1^{p-1}h_1\ot h_1^n$
mod $S^{2n-1}$ if $n>1$;
\item $d(v_1^\ell h_1^{n+1})\equiv-(\ell+n+1)v_1^\ell h_1\ot h_1^n$ mod $S^{2n-
1}$.
\end{enumerate}
\end{lem}

\section{New results about extensions}\label{extensions}
In this section, we show that $\vp_{2j-1}(X)$ is cyclic when $X$ is a sphere
bundle over a sphere with attaching map $\a_1$ or $\a_2$. This will be
crucial to our proof of Theorem \ref{main}. We also determine
the group structure of all groups $\vp_{2j-1}(X)$ when $X$ is an exceptional Lie
 group
for which the orders  $|\vp_*(X)|$ have been determined.

The first result of this section is the following, in which $B_k(2n+1,2n+kq+1)$
is an $S^{2n+1}$-bundle over $S^{2n+kq+1}$ with attaching map $\a_k$.
\begin{thm}\label{cyclic} Let $n>1$, and $k=1$ or $2$. Then
$v_{2j-1}(B_k(2n+1,2n+kq+1))$ and
$v_{2j}(B_k(2n+1,2n+kq+1))$ are isomorphic cyclic $p$-groups with exponent
$$\begin{cases}\min(n,2+\nu(j-n))&\text{if $j\equiv n$ mod $p(p-1)$}\\
\min(n+k(p-1),2+\nu(j-n-k(p-1)))&\text{if $j\equiv n$ mod $(p-1)$}\\
&\text{\qquad and $j\not\equiv n$ mod $p(p-1)$}\\
0&\text{otherwise}.\end{cases}$$
\end{thm}
\begin{pf} Let $B=B_k(2n+1,2n+kq+1)$. The determination of $v_{2j}(B)$ when
$k=1$ was made in \cite[1.3(2)]{BDMi}. In \cite[p.301]{F4}, $v_{2j}(B)$ was
determined when $k=2$, $n=4$, and $p=3$. The argument there adapts to the
general case in a straightforward fashion.

The cyclicity of $v_{2j-1}(B)$ when $k=1$ is proved similarly to
\cite[p.613]{DY}. It is easy when $|v_{2j}(S^{2n+q+1})|=p^{n+p-1}$, for then
$\partial:v_{2j}(S^{2n+q+1})\to v_{2j-1}(\sn)$ is surjective, and so
$v_{2j-1}(B)\approx v_{2j-1}(S^{2n+q+1})$ is cyclic.

Now consider the case when
$|v_{2j}(S^{2n+q+1})|<p^{n+p-1}$. By Theorem \ref{sph}(4), the class
$d(p^{m-n-p}h^m)\iota_{2n+q+1}$ has order $p$ in $E_2^2(S^{2n+q+1})$.
Here $m$ is an integer related to the stem of the class under consideration.
Since $\partial$ annihilates this class, there is $w\in VC^2(\sn)$ such that
$z:=d(p^{m-n-p}h^m)\iota_{2n+q+1}-w$ is a cycle in $E_2^2(B)$. We wish to show
that $pz$ is the image of a generator of $E_2^2(\sn)$.

We use the formula
\begin{equation}\label{bdryfor}
d(h^t\iota)=d(h^t)\iota+h^td(\iota),\end{equation} which was explained
as \cite[5.4]{DY}. This implies that
$$pz=d(p^{m-n-p+1}h^m\iota_{2n+q+1})-p^{m-n-p+1}h^m\ot h\iota_{2n+1}-pw.$$
The first term is a boundary. Using Lemma \ref{std}, the second term is,
mod terms that desuspend below $\sn$, $-v_1^{m-n-p+1}h^{n+p-1}\ot h\io_{2n+1}$,
and by Lemma \ref{yang}(2), this is, mod lower terms, $v_1^{m-n-p+1}h\ot
h^n\io_{2n+1}$, which is the leading term of a generator of
$E_2^2(\sn)$, by Theorem \ref{sph} or Proposition \ref{sph2}. Also,
as we shall show in the next paragraph, $pw$ desuspends to $S^{2n-1}$. Since
the double suspension from $E_2^2(S^{2n-1})$ to $E_2^2(\sn)$ is not surjective,
this implies that $pz$ is the image of a generator of $E_2^2(\sn)$.

One way to see that $w$ can be chosen so that $pw$ desuspends to $S^{2n-1}$ is
to note that $pz=d(p^{m-n-p+1}h^m)\iota_{2n+q-1}-pw$ is a cycle in
$E_2^2(B(2n-1,2n+q-1))$; i.e., multiplying by $p$ allows you to double
desuspend the whole equation.

The argument when $k=2$ is very similar. We will have $w\in VC^2(\sn)$
satisfying that
$z:=d(p^{m-n-2p+1}h^m)\iota_{2n+2q+1}-w$ is a cycle in $E_2^2(B)$,
and, as in the previous paragraph, $w$ can be chosen so that $pw$ double
desuspends. We obtain
$$pz=d(p^{m-n-2p+2}h^m\io_{2n+2q+1})-p^{m-n-2p+2}h^m\ot\a_2\io_{2n+1}-pw.$$
The first term is a boundary, the last term desuspends, while the middle term
is, mod terms that desuspend,
$v_1^{m-n-2p+2}h^{n+2p-2}\ot(2vh-ph^2)\io_{2n+1}$. The term with $ph^2$
desuspends, while the first term is, by \ref{yang}(2),
$2v_1^{m-n-p+1}h^{n+p-1}\ot h\io_{2n+1}$, which is the leading term of a
generator of $E_2^2(\sn)$. The class $\{pz\}$ generates $E_2^2(\sn)$ by
Proposition \ref{sph2}(1).
\end{pf}

In \cite{BDMi} and \cite{F4}, sphere bundles $X$ of the type covered by Theorem
\ref{cyclic} occurred as factors in product decompositions of exceptional Lie
groups (localized at a prime $p$). In those papers, we merely asserted the order
of the groups $\vp_{2j-1}(X)$, but we can now declare that they are cyclic.
There are a few other cases of factors $Y$ of exceptional Lie groups for which
only the order but not the group structure of $\vp_{2j-1}(Y)$ was given in
\cite{BDMi}, but we can now complete the determination of the $v_1$-periodic
homotopy of all torsion-free exceptional Lie groups by giving the group
structure in these cases. The following result handles all of these, and those
left unresolved in \cite{F4}.
\begin{prop} \label{complete} \begin{enumerate}
\item The $3$-primary groups $v_{2j-1}(B(11,15))$ and $v_{2j-1}(E_6/F_4)$,
which occur in \cite[1.2]{F4}, are cyclic.
\item The factors $B(2n+1,2n+q+1)$
which occur in $G_2$ for $p=5$, $F_4$ and $E_6$
for $5\le p\le 11$, $E_7$ for $11\le p\le 17$ and for $p=5$, and
$E_8$ for $11\le p\le 29$, as listed in \cite [1.1]{BDMi}, have $v_{2j-1}(B)$
cyclic of order given in \cite[1.3(2)]{BDMi}.
\item The spaces $B(11,23,35)$ and $B(23,35,47,59)$, which occur as factors
in $7$-primary $E_7$ and $E_8$, respectively, have $v_{2j-1}(B)$ cyclic of order
given in \cite[1.4]{BDMi}.
\item The spaces $B(3,11,19,27,35)$, $B(3,15,27)$, and $B(3,15,27,39)$, which
occur as factors of $5$-primary $E_7$, $7$-primary $E_7$, and $7$-primary $E_8$,
respectively, have $v_{2j-1}(B)\approx\bz/p\oplus\bz/p^{e-1}$, where
$e$ is the number given in \cite[1.4]{BDMi}.
\end{enumerate}\end{prop}
\begin{pf} The first two parts are immediate from Theorem \ref{cyclic}.
The first space in part 3 is a factor of $SU(18)$, and in the notation of
\cite[1.5]{DY} it has $N=5$ and $i=2$.
By \cite[1.9]{DY}, its
$v_{2j-1}(-)$-groups are cyclic. Similarly, the second space in part 3 is a
quotient of a factor $B=B(11,23,35,47,59)$ of $SU(30)$.  This factor $B$
has $N=5$ and $i=4$ in the notation of \cite[1.5]{DY}, and hence its groups
$v_{2j-1}(-)$ are cyclic by \cite[1.9]{DY}. Thus so are the groups of the
desired space $B(23,35,47,59)$, since $v_{2j-2}(S^{11})=0$ for values of $j$
under consideration.

The spaces in part 4 are also factors of $SU(n)$ and hence are covered by
\cite[1.9]{DY}. In the notation of
\cite[1.5]{DY}, these three spaces each have $N=1$, while $i=4$, 2, and 3,
respectively, and $\hat m>0$. Thus their
$v_{2j-1}(-)$ has a split $\bz/p$ by \cite[1.9]{DY}.
\end{pf}

\section{Discussion of $E_7/F_4$}\label{E7F4}
In this section we sketch a natural
approach to Theorem \ref{main}. Although we will not follow it exactly,
it is helpful in understanding the approach which we do employ.
Also, the result here about the convergence of the $v_1$-periodic UNSS for
$E_7/F_4$ will play a key role in our later deduction of $\vp_*(E_7)$.
Throughout the remainder of the paper, we will have $p=3$.

The fibration
\begin{equation}\label{EFfib}F_4\to E_7\to E_7/F_4\end{equation}
induces a long exact sequence of
$v_1$-periodic homotopy groups. The groups $v_*(F_4)$ were computed in
\cite{F4}, while $v_*(E_7/F_4)$ could be computed by
the methods of this paper.
Then we would need to determine the boundary homomorphism and
extensions in the exact sequence associated to (\ref{EFfib}).
This determination is complicated by the fact that the Bockstein
$\b$ is nonzero in
$H^*(F_4;\bz/3)$, which causes $BP_*(F_4)$ to be not a free $BP_*$-module, and
therefore the UNSS of $F_4$ cannot be calculated directly by known methods.
(In \cite{F4}, $v_*(F_4)$ was determined by a combination of topological and
UNSS methods.) Moreover, applying $\Omega$ to the fibration does not help much,
because $BP_*(\Om E_7)$ is not a free commutative algebra, and so we cannot
apply Theorem \ref{Vthm} to compute its $v_1$-periodic UNSS.
Hence UNSS methods cannot be used directly to analyze the exact sequence
in $v_*(-)$ associated to (\ref{EFfib}).

Our proof could be expedited slightly if we were assured of the validity of the
 following conjecture, due to Mimura.

\begin{conj} \label{Qconj} Localized at $3$, $E_7/F_4$ is spherically resolved
by spheres of dimension $19$, $27$, and $35$, and attaching maps $\a_2$. That
is, there is a fibration $S^{19}\to E_7/F_4\to B_2(27,35)$ and a fibration
$S^{27}\to B_2(27,35)\to S^{35}$, with attaching maps from $19$ to $27$ and
from $27$ to $35$ both equal to the element $\a_2$ which generates
$\pi_7(S^0)_{(3)}\approx\bz/3$.  \end{conj}

Although we cannot use this proposed topological description of $E_7/F_4$,
we can say enough about this space
to compute its $v_1$-periodic UNSS and prove that it converges to
$v_*(E_7/F_4)$. However, the specific results of this computation will not
be needed for the reasons cited earlier in this section, and the methods
will be applied again in computing the $v_1$-periodic UNSS of the space $Y_7$,
which will be our
approach to $v_*(E_7)$, and so we shall wait until the next section to use
them.

The following first steps toward proving Conjecture \ref{Qconj} will be useful
to us later. They were pointed out by Mimura.

\begin{prop}\label{a2h} $($a.$)$ $H^*(E_7/F_4;\bz)$ is an exterior algebra on
classes of dimension $19$, $27$, and $35$.
$($b.$)$ The $35$-skeleton of $E_7/F_4$ is $S^{19}\cup_{\pm\a_2}e^{27}\cup
_{\a_2}e^{35}$, where $\a_2$ generates $\pi_{7}(S^0)_{(3)}\approx\bz/3$.
\end{prop}

The proof of this proposition requires the following result of Kono and
Mimura.
\begin{prop}$($\cite{KM}$)$ \label{H^*E7}There is an algebra isomorphism
$$H^*(E_7;\bz_3)\approx\bz_3[e_8]/(e_8^3)\ot\Lambda[e_3,e_7,e_{11},e_{15},
e_{19},e_{27},e_{35}]$$
with only nonzero action of $\b$ or $\P^{p^r}$ on generators given by
$\b e_7=e_8$, $\b e_{15}=-e_8^2$, $\P^1e_3=e_7$, $\P^1e_{11}=e_{15}$,
$\P^1e_{15}=\pm e_{19}$, $\P^3e_7=e_{19}$, $\P^3e_{15}=e_{27}$.
\end{prop}

\begin{pf*}{Proof of Proposition \ref{a2h}} Part (a) was proved in
\cite[9.4]{MNT}. To prove part (b),
let $\Phi$ denote the secondary cohomology operation associated with
the relation $\P^1\b\P^1-\b\P^2-\P^2\b=0$. This secondary operation detects the
map $\a_2$ and satisfies $\P^3=\P^1\Phi$. (See \cite[p.353]{KM}.) In
\cite[7.2]{KM}, it is shown that $\Phi(\tilde e_{27})=\tilde e_{35}$ in
$H^*(\tilde E_7)$, where $\tilde E_7$ denotes the fiber of $E_7\to K(\bz,3)$,
from which it follows that $\Phi(e_{27})=e_{35}$ in $H^*(E_7)$.

Since $\P^1e_{15}=\pm e_{19}$ and $\P^3e_{15}=e_{27}$ in $H^*(E_7)$, we
can use a dual relation $\P^3=\Phi\P^1$ to deduce that
$\Phi(e_{19})=\pm e_{27}$. The dual relation is deduced by applying the
original relation in the $S$-dual, and then noting that $\P^1$, $\P^3$,
and $\Phi$ are all self-dual. Here duality is given by the
antiautomorphism of the Steenrod algebra, while $\Phi$ is self-dual
since it is defined by a symmetric Adem relation involving self-dual
terms.\end{pf*}

We close this section by proving the following result, which will be crucial
for us, since we will use it later to deduce that the
$v_1$-periodic UNSS of $Y_7$ converges to $v_*(Y_7)$.
\begin{thm}\label{convEF} The $v_1$-periodic UNSS of $E_7/F_4$ converges
to $v_*(E_7/F_4)$. Indeed, Theorem \ref{perUNSS} holds if $X=E_7/F_4$.
\end{thm}

Note that this result would be immediate from \ref{perUNSS}
if we knew that Conjecture
\ref{Qconj} were true. Instead, we must call upon the following result,
which was proved by Bendersky and Thompson at the request of the author.
\begin{thm}$($\cite{BT}$)$ \label{BTthm}
Suppose $X$ is a $K/p_*$-durable space with $K^*(X;\widehat
\bz_p)$ isomorphic as a
$\bz/2$-graded $p$-adic $\lambda$-ring to $\widehat\Lambda(M)$,
where $M=M_n$ is a
$p$-adic Adams module
which admits a sequence of epimorphisms of $p$-adic Adams modules
$$M_n@>p_n>> M_{n-1}@>p_{n-1}>> \cdots@>p_2>> M_1=M(2m_1+1)$$
with $\ker(p_i)=M(2m_i+1)$ for $2\le i\le n$. Here $M(2m+1)\approx
K^*(S^{2m+1};\widehat \bz_p)$ as a  $p$-adic Adams module,
and $\widehat\Lambda(M)$ denotes
the exterior algebra on $M$. Then the $(BP$-based$)$
$v_1$-periodic UNSS of $X$ converges to $v_*(X)$.\end{thm}

Actually, what is proved in \cite{BT} is that $X\to X\widehat{}$ induces
an isomorphism in $v_*(-)$. The target space denotes the $K/p$-completion.
In \cite{BT1}, it is proved that the $v_1$-periodic UNSS converges to
$v_*(X\widehat{})$, which then implies Theorem \ref{BTthm}.
The proof of this result in \cite{BT} relies heavily on the
work of Bousfield (\cite{Bo}). Bousfield defines a space $X$ to be
$K/p_*$-durable when its $K/p_*$-localization map induces an isomorphism in
$v_*(-)$. Theorem \ref{convEF} is an immediate consequence of Theorem
\ref{BTthm} and the following two results.
\begin{thm} There is an isomorphism of $\bz/2$-graded $p$-adic $\lambda$-rings
$$K^*(E_7/F_4;\widehat \bz_p)\approx\widehat\Lambda(M_3),$$
with short exact sequences of $p$-adic Adams modules
$$0\to M(35)\to M_3\to M_2\to 0\text{\quad and\quad}0\to M(27)\to M_2\to
M(19)\to 0.$$ \label{Mthm}
\end{thm}
\begin{thm}\label{durble} $E_7/F_4$ is $K/3_*$-durable.\end{thm}
\begin{pf*}{Proof of Theorem \ref{Mthm}} We use Proposition \ref{a2h}
to give the $E_2$-term of the Atiyah-Hirzebruch spectral sequence converging
to $K^*(E_7/F_4;\widehat \bz_p)$
as $\Lambda[x_{19},x_{27},x_{35}]\ot K\widehat \bz_p^*$.
The spectral sequence collapses to yield the claimed exterior algebra as
$K^*(E_7/F_4;\widehat \bz_p)$.
This collapsing can be deduced from Yagita's result
(\cite{Yag}) that there is a 3-local isomorphism
$$BP^*(E_7)\approx BP^*(F_4)\ot\Lambda[19,27,35],$$
or from Snaith's result (\cite{Sn}) that the spectral sequence
$$\on{Tor}_{R(G)}(\bz,R(H))\Rightarrow K^*(G/H)$$
collapses.

The claim about the decomposition of $M_3$ as a $p$-adic Adams module will
follow once we show that the generators of the exterior algebra
$K^1(E_7/F_4)$ satisfy $\psi^k(x_{35})=k^{17}x_{35}$,
$\psi^k(x_{27})=k^{13}x_{27}+\a_1x_{35}$, and $\psi^k(x_{19})=k^9x_{19}
+\a_2x_{27}+\a_3x_{35}$ for some integers $\a_1$, $\a_2$, and $\a_3$.
Note that $K^1(E_7/F_4)$ is spanned by $x_{19}$, $x_{27}$, $x_{35}$, and
$x_{19}x_{27}x_{35}$. We will show that the top cell of $E_7/F_4$, which
corresponds to this product class, splits off stably, and so cannot be involved
in Adams operations on the lower classes. Then the formula for the Adams
operations follows from the inclusions $S^{19}\to E_7/F_4$, $S^{27}\to
(E_7/F_4)/S^{19}$, and $S^{35}\to (E_7/F_4)/(E_7/F_4)^{(27)}$.

To prove the stable splitting, we argue similarly to \cite[1.1]{DMG}.
By \cite[3.3]{At}, the $S$-dual of the manifold $E_7/F_4$ is the Thom spectrum
of its stable normal bundle. However, $\widetilde{KO}(E_7/F_4)_{(3)}=0$, since
$E_7/F_4$ has no cells whose dimension is a multiple of 4. Thus the bottom
class splits off the Thom spectrum of the stable normal bundle, and dually the
top cell stably splits off the manifold itself.
\end{pf*}

The following proof is due to Pete Bousfield.
\begin{pf*}{Proof of Theorem \ref{durble}} In \cite{Bo}, Bousfield utilizes a
functor $\Phi$ from spaces to spectra, which he had introduced in earlier
papers. A map $f$ induces an isomorphism in $v_*(-)$ if and only if
$\Phi(f)$ is an equivalence. Let $X=E_7/F_4$, and consider the commutative
diagram
$$\begin{CD} \Phi(F_4)@>>>\Phi(E_7)@>>>\Phi(X)  \\
@VVV @VVV @VVV\\
\Phi((F_4)_{K/p})@>>>\Phi((E_7)_{K/p})@>>>\Phi(X_{K/p})
\end{CD}$$
Since $\Phi$ preserves fibrations, the top row is a fibration, and since
\cite[7.8]{Bo} states that $H$-spaces are $K/p_*$-durable, the first two
vertical arrows are equivalences. We will be done by the 5-lemma once we show
that the bottom row is a fibration.

By \cite[6.3]{Bo}, $K^*(G;\widehat
\bz_p)\approx\widehat\Lambda(P_G)$, where
$P_G=PK^1(G;\widehat \bz_p))$,
for $G=F_4$ or $E_7$, and by \cite[8.1]{Bo} $\Phi(G_{K/p})\simeq\Phi G$ is a
$K\widehat \bz_p^*$-Moore
spectrum $\cm(P_G/\psi^p,1)$, where $P_G/\psi^p$ is the quotient by the
injective action of the Adams operation. Similarly, by Theorem \ref{Mthm},
$$K^*(X_{K/p};\widehat \bz_p)\approx K^*(X;\widehat \bz_p)
\approx\widehat\Lambda(M_3),$$
and, since $X_{K/p}$, being $K/p$-local, is certainly $K/p_*$-durable,
we can apply \cite[8.1]{Bo} to obtain $\Phi(X_{K/p})\simeq \cm(M_3/\psi^p)$.
There is a short exact sequence of Adams modules, (e.g. from \cite{Yag})
$$0\to M_3\to PK^1(E_7)\to PK^1(F_4)\to 0$$
and hence a fiber sequence
$$\cm(P_{F_4}/\psi^p,1)\to \cm(P_{E_7}/\psi^p,1)\to \cm(M_3/\psi^p,1)$$
which is the bottom row of the commutative diagram considered above, showing
that it is a fibration, as desired.
\end{pf*}

\section{$E_2$ of periodic UNSS of $\Om E_7/Sp(2)$}\label{Y7}

In this long section, we calculate the periodic UNSS
 of $Y_7:=\Om E_7/Sp(2)$.
In Section \ref{E7}, we perform the  transition
from these results to $\vp_*(E_7)$.

We begin by recalling the following result of Harper, which we used in
\cite{F4}.
\begin{prop}$($\cite[4.4.1]{Har}$)$\label{HarK} There is a $3$-equivalence
$$F_4\approx K\times B(11,15),$$
where $K$ is a finite mod $3$ $H$-space satisfying
$$H^*(K;F_3)=\Lambda(x_3,x_7)\ot F_3[x_8]/(x_8^3),$$
with $x_7=\P^1x_3$ and $x_8=\b x_7$. Also, $B(11,15)$ is an $S^{11}$-bundle
over $S^{15}$ with $\P^1x_{11}=x_{15}$.
Moreover, there is a fibration $B(3,7)\to K\to W$, where $W$ is the Cayley
plane, and a fibration $S^7\to\Om W\to\Om S^{23}$.
\end{prop}

Because of the torsion in $H^*(K;\bz)$, and hence in $H^*(E_7;\bz)$,
we will work with loop spaces, and use the following result of
Hamanaka and Hara (\cite{HH}).
\begin{prop}\label{HaHa} The mod $3$ homology as Hopf algebras over $A$
satisfies
\begin{eqnarray*}
H_*(\Om F_4)&\approx& F_3[t_2,t_6,t_{10},t_{14},t_{22}]/(t_2^3)\\
H_*(\Om E_7)&\approx& F_3[t_2,t_6,t_{10},t_{14},t_{18},t_{22},t_{26},t_{34}]
/(t_2^3),
\end{eqnarray*}
with the only nonzero reduced coproducts being
$$\phibar(t_6)=-t_2^2\ot t_2-t_2\ot t_2^2$$
and \begin{eqnarray*}
\phibar(t_{18})&=&t_2^2t_6^2\ot t_2+t_2t_6^2\ot t_2^2-t_6^2\ot t_6
-t_2^2t_6\ot t_2t_6\\
&&\qquad -t_2t_6\ot t_2^2t_6-t_6\ot t_6^2+t_2^2\ot t_2t_6^2+t_2\ot t_2^2t_6^2.
\end{eqnarray*}
The only nonzero action of dual Steenrod operations $\P^{3^r}_*$ are
$\P^1_*(t_6)=t_2$, $\P^1_*(t_{14})=t_{10}$, $\P^1_*(t_{18})=\eps t_{14}-
t_2t_6^2$, $\P^1_*(t_{22})=\kappa t_6^3$, $\P^1_*(t_{26})=\eps t_{22}$,
$\P^1_*(t_{34})=-\eps t_{10}^3$, $\P^3_*(t_{18})=t_6$, $\P^3_*(t_{26})=t_{14}$,
and $\P^3_*(t_{34})=t_{22}$. Here $\eps=\pm1$ and $\kappa=\pm1$.
\end{prop}

Because of the relation $t_2^3=0$ in $H_*(\Om E_7)$, Theorem \ref{Vthm}
does not apply to $X=\Om E_7$. Instead, we will work with the space
$Y_7$ defined in the following theorem. We begin by noting (see
\cite[p.296]{F4}) that the space $B(3,7)$ which occurs in
\ref{HarK} is 3-equivalent to $Sp(2)$.
\begin{thm} \label{X7} Let $E_7/Sp(2)$
denote the quotient of the group inclusion
$Sp(2)\to F_4\to E_7$, and let $Y_7=\Om E_7/Sp(2)$. Then
$$H_*(Y_7;G)\approx
\Lambda[x_7]\ot G[x_{10},x_{14},x_{18},x_{22},x_{26},x_{34}]$$
for $G=\bz/3$ or $\bz_{(3)}$.
\end{thm}
\begin{pf} There is a commutative diagram of fibrations
\begin{eqnarray}
Y_7\quad&\to B(3,7)\to&E_7\nonumber\\
\da\quad&\da&\mapdown{=}\label{ffibr}\\
\Om E_7/F_4&\to\ F_4\ \to&E_7,\nonumber
\end{eqnarray}
 and this, together with the fibration
$B(3,7)\to F_4\to W\times B(11,15)$, which is a consequence of \ref{HarK},
implies there is a fibration
\begin{equation}
\Om W\times \Om B(11,15)\to Y_7\to \Om E_7/F_4.\label{X7fib}\end{equation}
The last fibration in \ref{HarK} determines $H_*(\Om W)$, and the Serre
spectral sequence of (\ref{X7fib}) collapses, yielding
the claim of the theorem. The collapsing is proved by observing that the only
possible differential on one of the three polynomial generators is
$d_{17}(x_{18})=\eps x_7\ot x_{10}$, but this has $\eps=0$ by consideration of
the map from (\ref{X7fib}) to the fibration
$$\Om W\times\Om B(11,15)\to B(3,7)\to F_4.$$
\end{pf}

We easily obtain the following consequence.
\begin{cor}\label{BPY7} $BP_*(Y_7)$ is a free commutative algebra on
 classes $x_7$, $x_{10}$, $x_{14}$,
$x_{18}$, $x_{22}$, $x_{26}$, and $x_{34}$, with $x_i\in BP_i(Y_7)$.
\end{cor}
\begin{pf} By \cite[12.1]{Ark}, the rationalization of $Y_7$
is homotopy equivalent to $K(Q,7)\times K(Q,10)\times\cdots\times K(Q,34)$.
Any differentials in the Atiyah-Hirzebruch spectral sequence
$$\Lambda[x_7]\ot
\bz_{(3)}[x_{10},x_{14},x_{18},x_{22},x_{26},x_{34}]
\ot BP_*\Rightarrow BP_*(Y_7)$$
must be seen rationally, and hence must be zero. That $x_7^2=0$ is deduced from
the inclusion $S^7\to Y_7$.
\end{pf}

By Theorem \ref{Vthm}, the UNSS of $Y_7$ can be calculated as the homology of
the unstable cobar complex. This complex splits as the direct sum of the
unstable cobar complex for $S^7$ plus the even-dimensional complex.
That is, we have
\begin{equation}\label{split12}E_2^{s,t}(Y_7)\approx\begin{cases}
E_2^{s,t}(S^7)&\text{if $t$ is odd}\\
\ext_\u^{s,t}(BP_*\langle x_{10},x_{14},x_{18},x_{22},x_{26},x_{34}\rangle)&
\text{if $t$ is even}
\end{cases}
\end{equation}
Our
work in this section will go into computing
$$v_1^{-1}
\ext_\u^{s,t}(BP_*\langle
x_{10},x_{14},x_{18},x_{22},x_{26},x_{34}\rangle).$$
This is the $v_1$-periodic $\ext$ which forms the $E_2$-term of the
$v_1$-periodic UNSS of $Y_7$.
In Section \ref{E7}, we will use Theorem \ref{BTthm} to show that this spectral
sequence converges to $v_*(Y_7)$. Throughout the remainder of the paper,
{\bf $E_2$ and $\ext_\u$ will always refer to their $v_1$-periodic versions},
unless explicitly stated to the contrary.

To compute the homology of the unstable cobar complex of $Y_7$, we will utilize
exact sequences in $\ext_\u(-)$ induced by the injective extension sequences
\begin{equation}\label{seq1}A(26)\to A(26,34)\to A(34),\end{equation}
\begin{equation}\label{seq2}A(18)\to A(18,26,34)\to A(26,34),\end{equation}
and
\begin{equation}\label{seq3}A(10,14)\ot A(22)\to BP_{\text{ev}}
(Y_7)\to A(18,26,34).\end{equation}
Each of these $A(-)$ is the subquotient of $BP_*(Y_7)$ on the generators of the
indicated dimensions. Each has an induced $\G$-coaction. The sequence
(\ref{seq3}) is closely related to the fibration
$$F_4\to E_7\to E_7/F_4,$$
with $F_4\approx B(11,15)\times K$.

By \cite[4.3]{BCR},
each of these three injective extension sequences yields a long exact sequence
when ordinary (unlocalized)
$\ext_\u(Q(-))$ is applied, and these $\ext$-groups are the homology of
the associated unstable cobar complexes.
The $v_1$-periodic $E_2$-term is the direct limit of
a direct system of $v_1$-power morphisms, and these commute with the morphisms
in the exact sequences just described. Since the direct limit of exact sequences
 is
exact, we obtain that there is an exact sequence of $v_1$-periodic $E_2$-terms.
As observed after Theorem \ref{perUNSS}, we can still work with the unstable
cobar complex, as long as we restrict attention to $v_1$-periodic classes.
 We will
abbreviate $\ext_\u(Q(A(n_1,\cdots,n_k)))$ as $E_2(n_1,\cdots,n_k)$, and the
associated unstable cobar complex as $C(n_1,\cdots,n_k)$.

In order to analyze $\partial$ in the long exact Ext sequences,
we will need the following
crucial result about the $\G$-coaction.
\begin{prop}\label{coass} If $M$ is a $\G$-comodule which as a $BP_*$-module
is free on $x_{10}$, $x_{14}$, $x_{18}$, $x_{26}$, and $x_{34}$, and if
\begin{eqnarray*}
\psi(x_{34})&=&1\ot x_{34}+\a_2\ot x_{26}+T_1\ot x_{18}+T_2\ot x_{14}+T_3\ot
x_{10}\\
\psi(x_{26})&=&1\ot x_{26}+\a_2\ot x_{18}+T_4\ot x_{14}+T_5\ot x_{10}\\
\psi(x_{18})&=&1\ot x_{18}+\a_1\ot x_{14}+T_6\ot x_{10}\\
\psi(x_{14})&=&1\ot x_{14}+\a_1\ot x_{10},
\end{eqnarray*}
then $T_6=\frac12h^2$, and, mod terms that desuspend lower than does the
indicated term, $T_1\equiv\frac12v^2h^2$, $T_2\equiv -5vh^4$,
$T_3\equiv\tfrac14vh^5$,
$T_4\equiv h^3$, and $T_5\equiv\tfrac14vh^3$.
\end{prop}

This proposition will be applied when $M$ is a quotient of $BP_*(Y_7)$.
The $\a_2$-terms in $\psi(x_{26})$ and $\psi(x_{34})$ are present there by
Proposition \ref{a2h}, since $\a_2$ is the cycle which detects the homotopy
class $\a_2$. We will see after the proof
that our application of this proposition to computing the homology of the
unstable cobar complex
would not be affected if a unit coefficient were present on $\a_2$.
Similarly, the $\a_1$-terms in $\psi(x_{18})$ and $\psi(x_{14})$ are present
because of $\P^1_*(t_{18})$ and $\P^1_*(t_{14})$ in \ref{HaHa}, and the
homology application would not be affected if they were multiplied by a unit.

\begin{pf} We begin with the determination of $T_1$. Using also
that $\psi{x_{26}}=1\ot x_{26}+\a_2\ot x_{18}$ mod lower terms,
the coassociativity formula
$(\psi\ot1)\psi(x_{34})=(1\ot\psi)\psi(x_{34})$
implies that $\psi(\a_2)=\a_2\ot1+1\ot\a_2$
and $\psi(T_1)=T_1\ot 1+1\ot T_1+\a_2\ot\a_2$. Now, $\a_2$ is given in Lemma
\ref{a2}, and one can verify that it is primitive. Let $\psibar$ be the reduced
coproduct in $BP_*BP$, defined by $\psi(y)=y\ot 1+1\ot y+\psibar(y)$.
We use the condition that $\psibar(T_1)=\a_2\ot\a_2$ to find $T_1$.

First, using Lemma \ref{a2} and Lemma \ref{std}(1), we compute
\begin{eqnarray*}\a_2\ot\a_2&=&(2vh-3h^2)\ot(2vh-3h^2)\\
&=&4vh\ot vh-6h^2\ot vh-6vh\ot h^2+9h^2\ot h^2\\
&=&4v^2h\ot h-18vh^2\ot h+18h^3\ot h-6vh\ot h^2+9h^2\ot h^2.
\end{eqnarray*}

Now $T_1$ must be a combination of the following five terms, whose $\psibar$ are
listed.
\begin{eqnarray}
h^4&\mapsto&4h^3\ot h+6h^2\ot h^2+4h\ot h^3.\label{psiT}\\
vh^3&\mapsto&3vh^2\ot h+3vh\ot h^2+v\ot h^3-1\ot vh^3\nonumber\\
&=&3vh^2\ot h+3vh\ot h^2+3h\ot h^3.\nonumber\\
v^2h^2&\mapsto&2v^2h\ot h+v^2\ot h^2-1\ot v^2h^2\nonumber\\
&=&2v^2h\ot h+6vh\ot h^2-9h^2\ot h^2.\nonumber\\
v^3h&\mapsto&v^3\ot h-1\ot v^3h=9v^2h\ot h-27vh^2\ot h+27h^3\ot h.\nonumber\\
h_2&\mapsto&
4h^3\ot h+6h^2\ot h^2+3h\ot h^3-vh\ot h^2-vh^2\ot h.\nonumber
\end{eqnarray}

We solve a system of linear equations for the coefficients of these five terms,
to see what combination $T_1$ can have $\psibar(T_1)=\a_2\ot\a_2$, as required.
We find that the desired term $T_1$ is given by
\begin{eqnarray*} T_1&=&\tfrac92h^4-6vh^3+2v^2h^2+c_1(-3h^4+vh^3+3h_2)
+c_2(-\tfrac{27}4h^4+9vh^3-\tfrac92 v^2h^2+v^3h),
\end{eqnarray*}
with $c_1$ and $c_2$ in $\bz_{(3)}$.
Replacing $3h$ by $v-\eta(v)$ at several places, this simplifies to
$$T_1=\tfrac12v^2h^2+L,$$
where $L$ desuspends to $S^3$.

The other $T_i$'s are determined similarly.
Coassociativity implies
\begin{eqnarray}
\psibar(T_6)&=&\a_1\ot \a_1\nonumber\\
\psibar(T_4)&=&\a_2\ot\a_1\nonumber\\
\psibar(T_5)&=&\a_2\ot T_6+T_4\ot\a_1\label{T3}\\
\psibar(T_2)&=&\a_2\ot T_4+T_1\ot\a_1\nonumber\\
\psibar(T_3)&=&\a_2\ot T_5+T_1\ot T_6+T_2\ot\a_1.\nonumber
\end{eqnarray}
That $T_6$ must equal $\frac12h^2$ is easily determined (since $\a_1=-h$).
To determine $T_4$, we write $\a_2\ot\a_1=-2vh\ot h+3h^2\ot h$, and note that
$\psibar$ acts as follows:
\begin{eqnarray*}
h^3&\mapsto&3h^2\ot h+3h\ot h^2\\
vh^2&\mapsto&2vh\ot h+3h\ot h^2\\
v^2h&\mapsto&6vh\ot h-9h^2\ot h
\end{eqnarray*}
Solving a system of equations for the coefficients of $h^3$, $vh^2$, and $v^2h$
yields
\begin{equation}\label{T4}T_4=h^3-vh^2+c(3h^3-3vh^2+v^2h),\end{equation}
with $c\in\bz_{(3)}$. All terms except the first are defined on $S^5$, and so
$T_4$ is as claimed.

We must have $\psibar(T_5)=vh\ot h^2-\frac32h^2\ot h^2-h^3\ot h+vh^2\ot h$.
The terms of which $T_5$ is a linear combination are the same as those in
$T_4$, which were listed with their $\psibar(-)$ in (\ref{psiT}). Solving this
system of equations yields that the combination whose $\psibar(-)$ is that
required of $T_5$ can be
$$-\tfrac14h^4+\tfrac13vh^3+c_1(-h^4+\tfrac13vh^3+h_2)+c_2(-\tfrac{27}4h^4
+9vh^3-\tfrac92v^2h^2+v^3h),$$
for any $c_1$ and $c_2$ in $\bz_{(3)}$. However, fractions with 3 in the
denominator do not lie in $\bz_{(3)}$.
 The only way to prevent this is to specify that $c_1$ must be of the
form $-1+3k$, with $k\in\bz_{(3)}$.  This yields
$$T_5=\tfrac34h^4-h_2+k(-3h^4+vh^3+3h_2)+L,$$
where $L$ (the $c_2$-term) desuspends to $S^3$. In two places,
we replace $3h^4$ by $vh^3-h^3v$, yielding
\begin{equation}\label{T5}T_5=\tfrac14vh^3+L',\end{equation}
as desired. In our determination of $T_5$, we should also take into account
the homogeneous part of (\ref{T4}), as it contributes to the $T_4\ot\a_1$-term
of $\psibar(T_5)$. When the resulting equations are solved,
we obtain an additional homogeneous part of $T_5$, equal to
$$c'(\tfrac34h^4-vh^3+\tfrac12v^2h^2)=c'(-\tfrac34vh^3-\tfrac34h^3v+\tfrac12v^2
h^2),$$
which desuspends farther than the leading term of (\ref{T5}). Thus $T_5$ is as
claimed.

Similarly, by (\ref{T3}) we must have (mod homogeneous terms that will be
considered below)
\begin{eqnarray*}\psibar(T_2)&=&(2vh-3h^2)\ot(h^3-vh^2)-(\tfrac92h^4-
6vh^3+2v^2h^2)\ot h\\
&=&2vh\ot h^3-3h^2\ot h^3-2v^2h\ot h^2+9vh^2\ot h^2-9h^3\ot h^2\\
&&-\tfrac92h^4\ot h+6vh^3\ot h-2v^2h^2\ot h.
\end{eqnarray*}
The terms that can comprise $T_2$ are listed below, with their $\psibar$.
\begin{eqnarray*}
h^5&\mapsto&5h\ot h^4+10h^2\ot h^3+10h^3\ot h^2+5h^4\ot h\\
vh^4&\mapsto&3h\ot h^4+4vh\ot h^3+6vh^2\ot h^2+4vh^3\ot h\\
v^2h^3&\mapsto&-9h^2\ot h^3+6vh\ot h^3+3v^2h\ot h^2+3v^2h^2\ot h\\
v^3h^2&\mapsto&27h^3\ot h^2-27vh^2\ot h^2+9v^2h\ot h^2+2v^3h\ot h\\
v^4h&\mapsto&-81h^4\ot h+108vh^3\ot h-54v^2h^2\ot h+12v^3h\ot h\\
vh_2&\mapsto&3vh\ot h^3+6vh^2\ot h^2+4vh^3\ot h-v^2h\ot h^2-v^2h^2\ot h+3h_1\ot
h_2\\
h_1h_2&\mapsto&3h\ot h^4+9h^2\ot h^3+10h^3\ot h^2+4h^4\ot h-vh\ot h^3\\
&&-2vh^2\ot h^2-vh^3\ot h+h_1\ot h_2+h_2\ot h_1\\
v_2h_1&\mapsto&24h^4\ot h-35vh^3\ot h+18v^2h^2\ot h-4v^3h\ot h+3h_2\ot h_1
\end{eqnarray*}
We solve a system of equations to find the combination of these terms having
$\psibar$ as desired. We obtain
\begin{eqnarray*}T_2&=&-\tfrac9{10}h^5+\tfrac32vh^4-
\tfrac23v^2h^3+c_1(\tfrac{81}5h^5-27vh^4+18v^2h^3-6v^3h^2+v^4h)\\
&&+c_2(-\tfrac{12}5h^5+7vh^4-\tfrac{17}3v^2h^3+2v^3h^2+vh_2-3h_1h_2+v_2h_1).
\end{eqnarray*}
As in the previous case, in order to prevent 3 in a denominator, we choose
$c_2=-1+3c$. This yields
$$T_2=\tfrac32h^5-\tfrac{11}2vh^4+5v^2h^3-2v^3h^2-vh_2+3h_1h_2-v_2h_1$$
plus two homogeneous terms which are defined on $S^7$. The first two terms in
$T_2$ combine to $\tfrac12vh^4-\tfrac12h^4v-\tfrac{11}2vh^4$, and so, mod terms
that are defined on $S^7$, we have $T_2\equiv -5vh^4$, as claimed.
We have omitted here consideration of homogeneous parts of $T_4$ and $T_1$
already obtained. These yield additional homogeneous terms in $T_2$ which are,
in fact, defined on $S^5$.

Finally we apply a similar method to determine $T_3$. It is again a matter of
solving a system of linear equations for the coefficients of the monomials that
can comprise $T_3$. We list the terms involved for the convenience of the
reader, who can quite easily check that our claimed $T_3$ does indeed have the
required coproduct. The lead term of this $T_3$ will play an important role in
our subsequent calculations. Indeed, it caused the answer for $v_*(E_7)$ to
turn out differently than the author had anticipated.

Momentarily ignoring some homogeneous parts, $T_3$ must satisfy
\begin{eqnarray} \psibar(T_3)&=&(2vh-3h^2)\ot(\tfrac34h^4-h_2)+(\tfrac92h^4-
6vh^3+2v^2h^2)\ot\tfrac12h^2\nonumber\\
&&+(\tfrac32h^5-\tfrac{11}2vh^4+5v^2h^3-2v^3h^2-vh_2+3h_1h_2-v_2h_1)\ot(-h).
\label{psiT3}\\
&=&-\tfrac94h^2\ot h^4+\tfrac94h^4\ot h^2-\tfrac32h^5\ot h+\tfrac32vh\ot h^4
-3vh^3\ot h^2\nonumber\\
&&+\tfrac{11}2vh^4\ot h+v^2h^2\ot h^2-5v^2h^3\ot h+2v^3h^2\ot h+3h^2\ot h_2
\nonumber\\
&&-3hh_2\ot h-2vh\ot h_2+vh_2\ot h+v_2h\ot h
\end{eqnarray}
We list the terms that can comprise $T_3$ and their coproducts.
\begin{eqnarray*}
h^6&\mapsto&6h\ot h^5+15h^2\ot h^4+20h^3\ot h^3+15h^4\ot h^2+6h^5\ot h\\
vh^5&\mapsto&3h\ot h^5+5vh\ot h^4+10vh^2\ot h^3+10vh^3\ot h^2+5vh^4\ot h\\
v^2h^4&\mapsto&-9h^2\ot h^4+6vh\ot h^4+4v^2h\ot h^3+6v^2h^2\ot h^2+4v^2h^3\ot
h\\
v^3h^3&\mapsto&27h^3\ot h^3-27vh^2\ot h^3+9v^2h\ot h^3+3v^3h\ot h^2+3v^3h^2\ot
h\\
v^4h^2&\mapsto&-81h^4\ot h^2+108vh^3\ot h^2-54v^2h^2\ot h^2+12v^3h\ot
h^2+2v^4h\ot h\\
v^5h&\mapsto&243h^5\ot h-405vh^4\ot h+270v^2h^3\ot h-90v^3h^2\ot h+15v^4h\ot
h\\
v^2h_2&\mapsto&3v^2h\ot h^3+6v^2h^2\ot h^2+4v^2h^3\ot h-v^3h\ot h^2-v^3h^2\ot
h\\
&&-9h^2\ot h_2+6vh\ot h_2\\
vhh_2&\mapsto&3vh\ot h^4+9vh^2\ot h^3+10vh^3\ot h^2+4vh^4\ot h-v^2h\ot h^3\\
&&-2v^2h^2\ot h^2-v^2h^3\ot h+3h\ot hh_2+vh\ot h_2+vh_2\ot h\\
h^2h_2&\mapsto&3h\ot h^5+12h^2\ot h^4+19h^3\ot h^3+14h^4\ot h^2+4h^5\ot h-vh
\ot h^4\\
&&-3vh^2\ot h^3-3vh^3\ot h^2-vh^4\ot h+2h\ot hh_2+h^2\ot h_2+2hh_2\ot h+h_2
\ot h^2\\
v_2h^2&\mapsto&24h^4\ot h^2-35vh^3\ot h^2+18v^2h^2\ot h^2-4v^3h\ot h^2+3h_2\ot
h^2+2v_2h\ot h\\
vv_2h&\mapsto&-72h^5\ot h+129vh^4\ot h-89v^2h^3\ot h+30v^3h^2\ot h-4v^4h\ot h\\
&&-9hh_2\ot h+3vh_2\ot h+3v_2h\ot h
\end{eqnarray*}

The solution of the resulting system of linear equations is
\begin{eqnarray}
T_3&=&\tfrac34h^6-\tfrac12v^2h^4+\tfrac12v^3h^3-\tfrac12v^2h_2+vhh_2-\tfrac
32h^2h_2+\tfrac12v_2h^2\label{T32}\\
&&+c_1(-\tfrac{81}2h^6+81vh^5-\tfrac{135}2v^2h^4+30v^3h^3-
\tfrac{15}2v^4h^2+v^5h)\nonumber\\
&&+c_2(9h^6-\tfrac{45}2vh^5+21v^2h^4-\tfrac{59}6v^3h^3+2v^4h^2+\tfrac12v^2h_2
\nonumber\\
&&\qquad -3vhh_2+\tfrac92h^2h_2-\tfrac32v_2h^2+vv_2h).\nonumber
\end{eqnarray}
The first term is rewritten as $\frac14(vh^5-h^5v)$, in order to see it
with a unit coefficient. All other terms desuspend to $S^9$.

The terms $T_1$, $T_5$, and $T_2$
which appear in the equation (\ref{T3}) for $\psibar(T_3)$ which gave rise to
the system of equations which we just solved have homogeneous parts whose
coefficients we do not know. For example, $T_5$ includes a summand of
$c(-3h^4+vh^3+3h_2)$. Thus added on to the RHS of (\ref{psiT3}) must be
$\a_2\ot c(-3h^4+vh^3+3h^2)$ and 7 other homogeneous parts arising similarly.
For each of these we solve a system of equations similar to the one just
solved, but with the RHS equal to the appropriate homogeneous term.
These give homogeneous summands to $T_3$. All resulting terms desuspend to
$S^7$, and so may be ignored. We spare the reader the details.
\end{pf}

The terms $\a_2\ot x_{26}$, $\a_2\ot x_{18}$, $\a_1\ot x_{14}$, and
$\a_1\ot x_{10}$ appear in the hypothesis of Proposition \ref{coass} because of
attaching maps in $\Om E_7$.
One might think that care is required as to the coefficients ($\pm1$)
of the $\a_2$ and $\a_1$ in Proposition \ref{coass}. However, this is not the
case. For if the four terms listed at the beginning of this paragraph are
multiplied by units $u_1$, $u_2$, $u_3$, and $u_4$, respectively, then
the terms $T_1$ to $T_6$ which are determined in Proposition \ref{coass}
are multiplied by units $u_1u_2$, $u_1u_2u_3$, $u_1u_2u_3u_4$,
$u_2u_3$, $u_2u_3u_4$, and $u_3u_4$, respectively. This can be seen by
consideration of the first part of the proof of \ref{coass}. For example,
we would have $\psibar(T_1)=u_1\a_2\ot u_2\a_2$.

The terms
$\a_1$, $\a_2$, and $T_i$ in \ref{coass}
will be used in the proofs of the theorems throughout the remainder of this
section to determine boundary morphisms in
exact sequences, and in pulling back terms whose boundary is 0. If units
$u_i$ were present as we are discussing here, it will only have the effect of
multiplying boundaries and pullbacks by unit amounts. The point is that all
terms in a boundary will be multiplied by the same unit, so that cancellation
due to different units cannot take place. For example, suppose that a term
$h^{I_1}x_{34}$ pulled back to $h^{I_1}x_{34}+h^{I_2}x_{26}$ in the case
where all $u_i=1$. Then, with units $u_i$ present, $h^{I_1}x_{34}$
pulls back to $h^{I_1}x_{34}+u_1h^{I_2}x_{26}$, and the boundary sends this
to $h^{I_1}\ot u_1u_2T_1x_{18}+u_1h^{I_2}\ot u_2\a_2x_{18}$, which is just $u_1u
_2$ times what it would have been. These uniform units do not affect whether
terms are zero, and hence can be ignored.

Now we can compute $E_2^{s,2j}(Y_7)$, dividing into cases depending upon the
parity and mod 9 value of $j$.  These will be delineated in Theorems \ref{1,7},
\ref{4}, \ref{0}, \ref{5,8}, \ref{2}, and \ref{lastthm}. Note that the exact
sequences in $E_2$ induced by (\ref{seq1}), (\ref{seq2}), and (\ref{seq3}),
together with (\ref{sphgp}), imply that if $t$ is even, then
$E_2^{s,t}(Y_7)=0$ unless $s=1$ or 2.

The first case is as follows.
\begin{thm}\label{1,7} If $j$ is odd, and $j\equiv1$ or $7$ mod $9$, then
$$E_2^{1,2j}(Y_7)\approx E_2^{2,2j}(Y_7)\approx \bz/3\oplus\bz/3^{\min(8,\nu(j-
43)+5)}.$$
\end{thm}
\begin{pf} Let $j$ be as in the theorem, and
$\nu=\nu(j-7)$. Formally, we obtain the result by computing first the
exact sequence in $E_2(-)$ associated to (\ref{seq1}), then that associated to
(\ref{seq2}), and then that associated to (\ref{seq3}). We know from
(\ref{shift}) and \ref{sph} that
\begin{equation}\label{sphgp}E_2^{s,2j}(2m)\approx\begin{cases}
\bz/3^{\min(m,\nu(j-m)+1)}&\text{if $j\equiv m$ mod 2, and $s=1$ or $2$}\\
0,&\text{otherwise} \end{cases}
\end{equation}  and we know from \cite[2.4]{BDMi} how to
compute $E_2(10,14)$ from $E_2(10)$ and $E_2(14)$. These are the building
blocks, but the glue is the boundary morphisms in the exact sequences, and
computing these requires much care.

A convenient way to picture the calculations is by Diagram \ref{17}, which we
think of as resembling an Adams spectral sequence chart.

\begin{diag}\label{17}
\begin{center}\begin{picture}(460,130)
\def\mp{\multiput}
\mp(89,15)(80,0){5}{\begin{picture}(40,110)
\put(-5,0){$E_2^{2,2j}$}
\put(26,0){$E_2^{1,2j}$}
\mp(1,18)(0,40){3}{$\bullet$}
\mp(31,18)(0,40){3}{$\bullet$}
\put(-1,35){$2$}
\put(29,35){$2$}
\mp(4,18)(0,40){2}{\line(0,1){17}}
\mp(34,18)(0,40){2}{\line(0,1){17}}
\mp(4,42)(0,40){2}{\line(0,1){18}}
\mp(34,42)(0,40){2}{\line(0,1){18}}
\mp(-5,58)(30,0){2}{$\bullet$}
\mp(0,42)(30,0){2}{\line(0,1){18}}
\end{picture}}
\put(90,90){$2$}
\put(120,90){$2$}
\put(170,90){$4$}
\put(200,90){$4$}
\put(250,90){$5$}
\put(280,90){$5$}
\put(330,90){$6$}
\put(360,90){$6$}
\put(410,90){$7$}
\put(440,90){$7$}
\put(17,30){$34$}
\put(17,50){$26$}
\put(22,70){$18$}
\put(22,90){$14$}
\put(22,110){$10$}
\put(7,70){$22$}
\put(23,38){\line(0,1){12}}
\put(19,58){\line(0,1){12}}
\mp(27,58)(0,20){3}{\line(0,1){12}}
\put(95,0){$\nu=1$}
\put(175,0){$\nu=3$}
\put(255,0){$\nu=4$}
\put(335,0){$\nu=5$}
\put(415,0){$\nu\ge6$}
\put(121,34){\vector(-1,3){26}}
\mp(200,40)(0,18){2}{\vector(-1,2){27}}
\mp(280,40)(0,18){2}{\vector(-1,2){27}}
\put(280,56){\vector(-2,3){27}}
\put(363,74){\vector(-3,4){30}}
\mp(361,33)(80,0){2}{\vector(-1,2){28}}
\mp(362,53)(0,4){2}{\vector(-3,4){29}}
\mp(442,53)(0,3){2}{\vector(-3,4){29}}
\mp(443,77)(0,17){2}{\vector(-3,2){29}}
\end{picture}\end{center}\end{diag}

Each $\bullet$ represents a $\bz/3$, and each integer $e$ represents a
$\bz/3^e$. These groups correspond to $E_2(2m)$, where $2m$ is the integer
 indicated on the left side of the diagram. The vertical lines
 indicate nontrivial extensions (multiplication by 3).
These are true because of the
$\a_1$ and $\a_2$ attaching maps and Theorem \ref{cyclic}. The positioning
of the 22-class is due to (\ref{X7fib}), i.e., that it is split away from the
10-, 14-, and 18-classes.

For example, the diagram for the case $\nu=3$ means that if $\nu(j-7)=3$
the boundary morphisms in (\ref{seq1}) and (\ref{seq2}) are 0, yielding
$E_2^{1,2j}(18,26,34)\approx E_2^{2,2j}(18,26,34)\approx\bz/3^4$, while in
(\ref{seq3})
$$E_2^{1,2j}(18,26,34)\mapright{\partial}E_2^{2,2j}(22)\oplus
E_2^{2,2j}(10,14)\approx \bz/3\oplus\bz/3^5$$
is 0 into the first summand and has image of order $3^2$ in the second summand.
Using either the exact sequence (\ref{seq3}) or the diagram,
this implies that in this case $$E_2^{1,2j}(Y_7)\approx E_2^{2,2j}(Y_7)\approx
\bz/3\oplus\bz/3^7.$$

The case $\nu=2$ omitted from Diagram \ref{17} has groups of order $3^3$,
i.e. labeled \lq\lq 3,'' on the 14-cell, and otherwise has the same groups
as do the other values of $\nu$. If $(j-7)/18\equiv1$ mod 3, then it has a
differential like that in the case $\nu=1$, while if $(j-7)/18\equiv2$ mod 3,
then it has no nonzero differentials.

Of course, we still have to verify that the differentials are as claimed
in Diagram \ref{17} and the above paragraph describing the case $\nu=2$.
The reader can easily verify that this will imply Theorem \ref{1,7}.
The $\bz/3$ on the 22-cell splits for algebraic reasons.

We analyze the differentials by the methods used extensively in
\cite{DY} and \cite{exp}, involving the unstable cobar complex.
One convention is that we often omit writing
powers of $v_1$ on the left; they can always be determined by consideration of
total degree. The boundary $E_2^{1,2j}(34)\mapright{\partial}
E_2^{2,2j}(26)$ sends the
generator $h\io_{34}$ to $h\ot\a_2\io_{26}=h\ot(hv+vh)\io_{26}$. This is
obtained from $\psi(x_{34})=1\ot x_{34}+\a_2\ot x_{26}$ in \ref{coass},
and from \ref{a2}. The relationship of $\partial$ with the coaction is
standard; see, e.g., \cite[2.7]{exp}.

We use \ref{std} to write $h\ot vh=vh\ot h-3h^2\ot h$,
and $3h^2\ot h=h^2\ot(v-\eta v)=0$. Also, $h\ot hv$ is defined on $S^1$ and
hence is 0 in $E_2(-)$.
So the image of $\partial$ equals $v^{\text{pwr}}h\ot h\io_{26}$.
By \ref{sph2}(second part of (2)), this cycle equals $d((h^3+L)\io_{26})$,
with $\nu=\nu(j-13)=1$ in
\ref{sph2}. Here we have
omitted a unit coefficient, which will be done routinely unless the
coefficient plays a significant role. Here $L$ desuspends {\bf lower} than the
associated term, in this case $h^3$, a notation that will be employed
frequently, with the $L$'s sometimes adorned with primes to distinguish them
from one another.

Thus the generator $h\io_{34}$ pulls back to $h\io_{34}+(h^3+L)\io_{26}$
in $E_2^{1,2j}(26,34)$, and $\partial$ in (\ref{seq2}) sends this to
$(h\ot(\frac12v^2h^2+L')+(h^3+L)\ot\a_2)\io_{18}\in E_2^{2,2j}(18)$.
Here we have used \ref{coass}. The leading term here is
$$\tfrac12h\ot v^2h^2=\tfrac12(\eta v)^2h\ot h^2= \tfrac12(v-3h)^2h\ot h^2,$$
which has leading term $\frac12 h\ot h^2$. Note how $v$'s on the left are
absorbed into other unstated $v$'s. By \ref{sph2}(2), this equals
$d((h^3+L'')\io_{18})$ (omitting unit coefficients), and so our generator pulls
back to
$$z\equiv h\io_{34}+h^3\io_{26}+h^3\io_{18}\in E_2^{1,2j}(18,26,34).$$
Here, and subsequently, \lq\lq$\equiv$'' will mean \lq\lq mod $L$,'' with the
lower terms varying from term to term.

We analyze the two components of $\partial(z)$ in the exact sequence of $E_2$
derived from (\ref{seq3}).
We begin by showing that the component $\partial_2$
into $E_2^{2,2j}(22)$ is 0. We have
$\partial_2(z)=(h\ot (h^3+L)+h^3\ot h)\io_{22}$. Here we use
the $\a_1$ attaching map from 22 to 26 in $\Om E_7$,
which causes the $\ot h\io_{22}$. The $\ot(h^3+L)$ is obtained by the same
calculation that gave $T_4$ in Proposition \ref{coass}. But these terms don't
even matter very much, for such terms desuspend far below $S^{22}$, and hence
are 0 in $E_2^{2,2j}(22)\approx\bz/3$. Here we use a fact that we will use
frequently, essentially from \ref{sph}(5), that if $E_2^{2,2j}(2n+\eps)\approx
\bz/p$, then an element in it which is in the image of the double desuspension
is 0.

Similar, although much more delicate, considerations apply to obtaining
the other component
\begin{equation}\label{par1}
\partial_1:E_2^{1,2j}(18,26,34)\to E_2^{2,2j}(10,14).\end{equation}
First we determine the composite when $\partial_1$ is followed (by $\rho$)
into $E_2^{2,2j}(14)$.
 Using Proposition \ref{coass} and the
usual relationship between the coaction and the boundary morphism, we obtain
\begin{equation}\label{rhod}
\rho\partial_1(z)\equiv(h\ot vh^4+h^3\ot h^3+h^3\ot h)\io_{14}.
\end{equation}
Here all terms except $h\ot vh^4$ desuspend to $S^7$, while, mod $S^7$,
$h\ot vh^4\equiv h\ot h^4$. If $\nu\ge3$, then, by \ref{sph2}(1), this element
has order $3^{\min(4,\nu-2)}$ in $E_2^{2,2j}(14)$, and this is as claimed in
Diagram \ref{17}, with the arrows above the lowest one being a consequence
of the lowest one and the extensions.

If $\nu=1$, then by \ref{sph2}(2)
and (\ref{rhod}), $\rho\partial_1(z)= d((h^6+L)\io_{14})$.
Thus $z$ pulls back to
$$z'\equiv h\io_{34}+h^3\io_{26}+h^3\io_{18}+h^6\io_{14}\in C(14,18,26,34).$$
Using \ref{coass}, this satisfies
$$\partial(z')\equiv (h\ot vh^5+h^3\ot vh^3+h^3\ot h^2+h^6\ot h)\io_{10}.$$
Here there can be \lq\lq lower'' terms associated with the factor on either
side of the tensor sign, omitted $v$'s occur only on the left, and, as usual,
unit coefficients are omitted. All terms here except the first desuspend to
$S^9$, while that term generates $E_2^{2,2j}(10)$, so the image of
$\partial_1$ in this case has order $3$ in $E_2^{2,2j}(10,14)$, as claimed.

Finally we consider the delicate case when $\nu=2$. In this case, there are two
terms with the potential to cancel, and so we must keep track of unit
coefficients. We write $j=7+18c$, with $c\not\equiv0$ mod 3. As before,
$E_2^{1,2j}(18,26,34)$ is generated by
$z\equiv h\io_{34}+h^3\io_{26}+h^3\io_{18}$. The
unit coefficients of the second and third terms will not be important, and so
are
omitted. The leading term of $\rho\partial_1(z)$ in $E_2^{2,2j}(14)$ is,
by \ref{coass}, $h\ot(-5vh^4)\io_{14}\equiv h\ot h^4\io_{14}$, where we have
used that $-5\equiv1$ mod 3. By Lemma \ref{yang}(3), $d(h^7)\equiv
-9ch\ot h^6$ in
this stem. ($(\ell+n+1)$ of the lemma multiplied by $2(p-1)$ equals $2j-14$.)
Thus, since $9h^6\equiv h^4$, we obtain $\rho\partial_1(z)=d((-
\frac1ch^7+L)\io_{14})$,
and so $z$ pulls back to
$$z'\equiv h\io_{34}+h^3\io_{26}+h^3\io_{18}+\tfrac1ch^7\io_{14}.$$
This satisfies
$$\partial(z')\equiv (h\ot\tfrac14vh^5+uh^3\ot vh^3+u'h^3\ot h^2+\tfrac1ch^7\ot
(-h))\io_{10}$$
with $u$ and $u'$ units in $\bz_{(3)}$.
The middle terms desuspend, while the first and last combine,
using \ref{yang}(2), to give
$\frac14+\frac1c$ times the generator of $E_2^{2,2j}(10)\approx\bz/3$.
This is nonzero if $c\equiv1$ mod 3, and 0 if $c\equiv2$ mod 3, as claimed
in the paragraph earlier in the proof which described the case $\nu=2$.
\end{pf}

The statement and proof for the case $j\equiv4$ mod 9 are quite similar to the
cases just completed.
\begin{thm}\label{4} If $j$ is odd, and $j\equiv4$  mod $9$, then
$$E_2^{1,2j}(Y_7)\approx E_2^{2,2j}(Y_7)
\approx \bz/3\oplus\bz/3^{\min(14,\nu(j-13-4\cdot 3^8)+5)}.$$
\end{thm}
\begin{pf} Let $j$ be as in the theorem, and
$\nu=\nu(j-13)$. As in the previous theorem, the way in which the result stated
in the theorem is obtained is most conveniently expressed in a diagram.

\begin{diag}\label{d4}
\begin{center}\begin{picture}(460,130)
\def\mp{\multiput}
\mp(89,15)(80,0){5}{\begin{picture}(40,110)
\put(-5,0){$E_2^{2,2j}$}
\put(26,0){$E_2^{1,2j}$}
\mp(1,18)(0,40){3}{$\bullet$}
\mp(31,18)(0,40){3}{$\bullet$}
\put(1,75){$2$}
\put(31,75){$2$}
\mp(4,18)(0,40){2}{\line(0,1){17}}
\mp(34,18)(0,40){2}{\line(0,1){17}}
\mp(4,42)(0,40){2}{\line(0,1){18}}
\mp(34,42)(0,40){2}{\line(0,1){18}}
\mp(-5,58)(30,0){2}{$\bullet$}
\mp(0,42)(30,0){2}{\line(0,1){18}}
\end{picture}}
\put(78,50){\frame{$\nu+1$}}
\put(108,50){\frame{$\nu+1$}}
\put(166,50){$10$}
\put(196,50){$10$}
\put(248,50){$11$}
\put(278,50){$11$}
\put(326,50){$12$}
\put(356,50){$12$}
\put(406,50){$13$}
\put(436,50){$13$}
\put(17,30){$34$}
\put(17,50){$26$}
\put(22,70){$18$}
\put(22,90){$14$}
\put(22,110){$10$}
\put(7,70){$22$}
\put(23,38){\line(0,1){12}}
\put(19,58){\line(0,1){12}}
\mp(27,58)(0,20){3}{\line(0,1){12}}
\put(95,0){$\nu\le7$}
\put(175,0){$\nu=9$}
\put(255,0){$\nu=10$}
\put(335,0){$\nu=11$}
\put(415,0){$\nu\ge12$}
\put(121,34){\vector(-1,3){26}}
\mp(200,40)(0,18){2}{\vector(-1,2){27}}
\mp(280,39)(0,19){2}{\vector(-1,2){27}}
\put(280,55){\vector(-2,3){28}}
\mp(362,58)(80,0){2}{\vector(-1,2){27}}
\mp(362,52)(0,3){2}{\vector(-3,4){28}}
\mp(442,52)(0,3){2}{\vector(-3,4){28}}
\put(360,34){\vector(-3,4){30}}
\mp(441,37)(-3,17){2}{\vector(-3,2){27}}
\end{picture}\end{center}\end{diag}

The omitted case $\nu=8$ is like the case $\nu\le7$ if $(j-
13)/(2\cdot3^8)\equiv1$ mod 3, while it has all differentials 0 if
$(j-13)/(2\cdot3^8)\equiv2$ mod 3. In most cases, the $\bz/3$ from
the 22-class splits for algebraic reasons. The splitting in the cases
when $\nu\ge11$
require a bit of care, which will be dealt with later in the proof.
The boundary in the diagram in these cases is meant to be hitting the
sum of the classes on the 22 and the 18.

We begin with the case $\nu\le7$. We start as in the proof of \ref{1,7},
but this time the boundary of $h\io_{34}$ in $E_2^{2,2j}(26)$ is
$d((h^{\nu+2}+L')\io_{26})$, by \ref{sph2}(2). (In the proof of \ref{1,7}, we
had $\nu(j-13)=1$.) Thus
the generator of $E_2^{1,2j}(26,34)$ equals, mod lower terms,
$h\io_{34}+h^{\nu+2}\io_{26}$. The next term is found by writing
\begin{equation}\label{nu2}
(h\ot(\tfrac12v^2h^2+L)+(h^{\nu+2}+L')\ot(2vh-3h^2))\io_{18}
\end{equation} as a boundary in
the unstable cobar complex. The first term will dominate if $\nu\le4$, while
the second term will dominate if $4\le\nu\le7$. (If $\nu=4$, there could be
cancellation that would cause it to desuspend even lower, but that won't affect
the final result.) We obtain that
$$z\equiv h\io_{34}+h^{\nu+2}\io_{26}+h^{\max(3,\nu-1)}\io_{18}$$
generates $E_2^{1,2j}(18,26,34)$. The $h^{\nu-1}\io_{18}$ when $\nu\ge4$
is obtained since
$$h^{\nu+2}\ot vh\equiv h^{\nu}\ot h\equiv h\ot h^{\nu-2}=d(h^{\nu-1}),$$
using \ref{yang}(1,2,3).
Now $$\rho\partial_1(z)\equiv (h\ot vh^4+h^{\nu+2}\ot h^3
+h^{\max(3,\nu-1)}\ot h)\io_{14},$$
which has leading term a multiple, $k$, of $h\ot h^4$. By \ref{sph2}(2)
this is $d(kh^6\io_{14})$ since $\nu(2j-14)=1$, and so $z$ pulls back to
$z'\equiv z-kh^6\io_{14}$ in $C(14,18,26,34)$. The leading term of
$\partial(z')$ in $E_2^{2,2j}(10)$ is $h\ot vh^5$, which is a generator.
It is also important to know here that
$E_2^{1,2j}(18,26,34)\mapright{\partial} E_2^{2,2j}(22)$ is 0, for if it were
nonzero then $E_2^{2,2j}(Y_7)$ would be cyclic. The leading term of this
$\partial$ is $h^{\nu+2}\ot h\io_{22}$ which desuspends and hence is 0 in $E_2$.

Next we consider the case $\nu=8$, in which we have to keep track of
unit coefficients because of the possibility of two cancelling terms.
Let $j-13=2\cdot 3^8c$, with $c\not\equiv0$ mod 3. By \ref{yang}(3),
we have $d(h^{10}\io_{26})\equiv-3^8ch\ot h^9\io_{26}\equiv -ch\ot h\io_{26}$,
where the second step utilizes $3^8h^8=(v-\eta v)^8$. Then $E_2^{1,2j}(34)
\mapright{\partial}E_2^{2,2j}(26)$ sends the generator to
$h\ot(-h)\io_{26}=d((\frac1ch^{10}+L)\io_{26})$,
and so the generator pulls back to
$z\equiv h\io_{34}-\frac1ch^{10}\io_{26}$. The leading term of $\partial(z)$ in
$E_2^{2,2j}(18)$ is $-\frac1ch^{10}\ot2vh\io_{18}$, which, using
\ref{yang}(1,2), is equivalent to $-\frac2ch^8\ot h\io_{18}\equiv \frac2ch\ot
h^6\io_{18}$. By \ref{yang}(3), $$d(h^7\io_{18})\equiv
-\tfrac12(j-9)h\ot h^6\io_{18}
=-(2+3^8c)h\ot h^6\io_{18}\equiv -c\partial(z).$$
Thus $z$ pulls back to $z'\equiv z+\frac1ch^7\io_{18}$. The leading term of
$\rho\partial_1(z')$ in $E_2^{2,2j}(14)$ is $\frac1ch^7\ot(-h)\io_{14}$.
By \ref{yang}(3) again, $$d(h^7\io_{14})\equiv -\tfrac12(j-7)h\ot h^6\io_{14}
\equiv -h\ot 3h^6\io_{14}\equiv-h\ot h^5\io_{14}\equiv-c\partial(z'),$$
where we have used \ref{yang}(1) at the last step. Thus $z'$ pulls back to
$z''\equiv z'+\frac1ch^7\io_{14}$ in $C(14,18,26,34)$. There are two leading
terms in $\partial(z'')\in E_2^{2,2j}(10)$. These are $h\ot\frac14vh^5$
and $\frac1ch^7\ot(-h)$. They combine to give $\frac14+\frac1c$ times a
generator, and this is 0 if $c\equiv2$ mod 3, and nonzero if $c\equiv 1$ mod 3,
as claimed. The $\partial$ into the $22$-part is 0 as in the case
$\nu\le7$.

If $9\le\nu\le10$, the situation is much easier. Similarly to the previous
cases, but ignoring units, the generator of $E_2^{1,2j}(18,26,34)$ is
$z\equiv h\io_{34}+h^{\nu+2}\io_{26}+h^{\nu-1}\io_{18}$. The leading term of
$\rho\partial_1(z)$ in $E_2^{2,2j}(14)$ is $h^{\nu-1}\ot h\io_{14}$, which is
a generator if $\nu=10$, and is 3 times the generator if $\nu=9$. The boundary
$E_2^{1,2j}(18,26,34)\mapright{\partial_2}E_2^{2,2j}(22)$ is 0 because its
leading
term is $h^{\nu+2}\ot h\io_{22}$ which is 0 in $E_2$ for $\nu\le10$.

When $\nu=11$, the boundary from $E_2^{1,2j}(26,34)$ to $E_2^{2,2j}(18)$
is now nonzero. Indeed, its image, given in (\ref{nu2}), has leading term
$$h^{13}\ot 2vh\io_{18}\equiv 2h^{11}\ot h\io_{18}\equiv h\ot h^9\io_{18},$$
which is a generator. Here we have used \ref{yang}(1) and \ref{yang}(2).
The boundary from $E_2^{1,2j}(26,34)$ to $E_2^{2,2j}(22)$ is also nonzero
since the generator $z\equiv h\io_{34}+h^{13}\io_{26}$
satisfies $\partial(z)\equiv h^{13}\ot h
\io_{22}$, and this is a generator. The chart would
then suggest (accurately) that the boundary hits into the sum of the two
classes, and the extension is also into this sum. One way to formalize this
uses the exact sequence
\begin{equation}\label{R}
E_2^{1,2j}(26,34)\mapright{\partial}E_2^{2,2j}(10,14,18,22)\to E_2^{2,2j}(Y_7)
\to E_2^{2,2j}(26,34).\end{equation}
The first and last groups are $\bz/3^{13}$, while
the  second is $\bz/3\oplus\bz/3^4$. The boundary $\partial$ hits the sum of
the two generators. There is a cycle representative
$z$ in $E_2^{2,2j}(Y_7)$ which projects to
an element of order $3$ in $E_2^{2,2j}(26,34)$ and satisfies that
3 times this generator is the image of the sum of the two generators of
$E_2^{2,2j}(10,14,18,22)$. This implies $E_2^2(Y_7)
\approx \bz/3\oplus \bz/3^{13}$.

Actually, a little bit more care is required here with regard to coefficients
of the generators. It is conceivable that the boundary could hit the sum of
generators but the extension be into their difference, and then the extension
group would be cyclic of order $3^{14}$. What really happens is that, if $c$ is
defined as before by $j-13=2\cdot 3^{11}c$, then a
generator $z\equiv h\io_{34}-\frac1ch^{13}\io_{26}$ satisfies
\begin{equation}\label{form1}
\partial(z)\equiv \tfrac1ch^{13}\ot h\io_{22}-\tfrac1ch^{13}\ot 2vh\io_{18},
\end{equation}
while on the other hand, by the argument of Theorem \ref{cyclic}, the element
$d(3^{j-27}h^{j-13})\io_{26}$ of order 3 extends to a cycle $z'$
in $C(18,22,26)$ such that, mod classes that desuspend farther, $3z'$ is
homologous to
\begin{equation}\label{form2}
3^{j-26}h^{j-13}\ot(-h)\io_{22}+3^{j-26}h^{j-13}\ot2vh\io_{18}.\end{equation}
The classes in (\ref{form1}) and (\ref{form2}) are clearly unit multiples of
one another. In each case, we use $h^{13}\ot vh\equiv h^{11}\ot h$ to see that
the second term is a generator.

This completes the case $\nu=11$.  The case $\nu\ge12$ is very similar.
Actually it is a bit easier, for the consideration of the previous paragraph
need not be addressed, since the initial differential hits into a cyclic group.
\end{pf}

The case $j\equiv0$ mod 3 introduces no new ideas.
\begin{thm}\label{0} If $j$ is odd, and $j\equiv0$  mod $3$, then
$$E_2^{1,2j}(Y_7)\approx E_2^{2,2j}(Y_7)
\approx \bz/3\oplus\bz/3^{\min(10,\nu(j-9-2\cdot 3^5)+4)}.$$
\end{thm}
\begin{pf} Let $j$ be as in the theorem, and $\nu=\nu(j-9)$.
We will show that the diagram encapsulating the exact sequences of
(\ref{seq1}), (\ref{seq2}), and (\ref{seq3}) is as depicted in
Diagram \ref{charts2} for certain values of $\nu$. This diagram, together with
the subsequent discussion of what happens for values of $\nu$ not included in
the diagram, implies Theorem \ref{0}.

\begin{diag}\label{charts2}
\begin{center}
\begin{picture}(460,190)
\put(21,41){$34$}
\put(21,71){$26$}
\put(36,101){$18$}
\put(6,101){$22$}
\put(36,131){$14$}
\put(36,161){$10$}
\put(25,55){\line(0,1){10}}
\put(40,115){\line(0,1){10}}
\put(40,145){\line(0,1){10}}
\put(20,80){\line(-1,2){10}}
\put(30,80){\line(1,2){10}}
\put(97,44){$\bullet$}
\put(87,105){$\bullet$}
\put(107,100){$8$}
\put(107,163){$\bullet$}
\put(97,72){$\bullet$}
\put(107,132){$\bullet$}
\put(100,48){\line(0,1){22}}
\put(110,108){\line(0,1){22}}
\put(110,139){\line(0,1){24}}
\put(98,80){\line(-1,3){8}}
\put(102,80){\line(1,3){8}}
\put(147,44){$\bullet$}
\put(137,105){$\bullet$}
\put(157,100){$8$}
\put(157,163){$\bullet$}
\put(147,72){$\bullet$}
\put(157,132){$\bullet$}
\put(150,48){\line(0,1){22}}
\put(160,108){\line(0,1){22}}
\put(160,139){\line(0,1){24}}
\put(148,80){\line(-1,3){8}}
\put(152,80){\line(1,3){8}}
\put(277,44){$\bullet$}
\put(267,105){$\bullet$}
\put(287,100){$7$}
\put(287,163){$\bullet$}
\put(277,72){$\bullet$}
\put(287,132){$\bullet$}
\put(280,48){\line(0,1){22}}
\put(290,108){\line(0,1){22}}
\put(290,139){\line(0,1){24}}
\put(278,80){\line(-1,3){8}}
\put(282,80){\line(1,3){8}}
\put(229,44){$\bullet$}
\put(219,105){$\bullet$}
\put(239,100){$7$}
\put(239,163){$\bullet$}
\put(229,72){$\bullet$}
\put(239,132){$\bullet$}
\put(232,48){\line(0,1){22}}
\put(242,108){\line(0,1){22}}
\put(242,139){\line(0,1){24}}
\put(230,80){\line(-1,3){8}}
\put(234,80){\line(1,3){8}}
\put(407,44){$\bullet$}
\put(397,105){$\bullet$}
\put(407,100){\frame{$\nu+1$}}
\put(417,163){$\bullet$}
\put(407,72){$\bullet$}
\put(417,132){$\bullet$}
\put(410,48){\line(0,1){22}}
\put(420,108){\line(0,1){22}}
\put(420,139){\line(0,1){24}}
\put(408,80){\line(-1,3){8}}
\put(411,77){\line(1,3){8}}
\put(357,44){$\bullet$}
\put(347,105){$\bullet$}
\put(357,100){\frame{$\nu+1$}}
\put(367,163){$\bullet$}
\put(357,72){$\bullet$}
\put(367,132){$\bullet$}
\put(360,48){\line(0,1){22}}
\put(370,108){\line(0,1){22}}
\put(370,139){\line(0,1){24}}
\put(358,80){\line(-1,3){8}}
\put(361,77){\line(1,3){8}}
\put(146,46){\vector(-1,2){31}}
\put(146,74){\vector(-1,2){31}}
\put(157,100){\vector(-2,3){44}}
\put(276,45){\vector(-1,3){30}}
\put(276,75){\vector(-1,3){30}}
\put(405,45){\vector(-1,4){30}}
\put(93,25){$E_2^{2,2j}$}
\put(223,25){$E_2^{2,2j}$}
\put(353,25){$E_2^{2,2j}$}
\put(143,25){$E_2^{1,2j}$}
\put(273,25){$E_2^{1,2j}$}
\put(403,25){$E_2^{1,2j}$}
\put(116,5){$\nu=7$}
\put(246,5){$\nu=6$}
\put(376,5){$\nu\le4$}
\end{picture}
\end{center}
\end{diag}

If $\nu\ge8$, then the group corresponding to the 18-cell has order $3^9$,
and by \cite[2.4]{BDMi} there is a nonzero boundary morphism from
$E_2^{1,2j}(26)$ to $E_2^{2,2j}(18)$ (hitting the element of order $3$, of
course), and three other boundary morphisms (one below and two above it)
follow from it by the extensions in a diagram similar to that of \ref{charts2}.

We will show below that the case $\nu=5$ is like the case $\nu\le4$ if
$(j-9)/(2\cdot3^5)\equiv2$ mod 3, while it has no differentials if
$(j-9)/(2\cdot3^5)\equiv 1$ mod 3. But first we establish that the cases
in Diagram \ref{charts2} are as depicted.

Let $\nu=7$. The nonzero differential from $E_2^{1,2j}(34)\approx\bz/3$
to $E_2^{2,2j}(18)\approx\bz/3^8$
is established similarly to that in the case $\nu=11$ in the preceding
theorem. Indeed, the boundary $E_2^{1,2j}(34)\mapright{\partial}E_2^{2,2j}(26)$
sends the generator to
$$h\ot\a_2\io_{26}=h\ot(2vh-3h^2)\io_{26}=(2vh\ot h-6h^2\ot h-3h\ot
h^2)\io_{26}=d(vh^2-2h^3)\io_{26},$$
and so it pulls back to
$h\io_{34}+(2h^3-vh^2)\io_{26}$. This in turn has boundary
$(h\ot(\frac12v^2h^2+L)+(2h^3-vh^2)\ot\a_2)\io_{18}$, whose leading term
$\frac12h\ot v^2h^2\io_{18}\equiv \frac12h\ot h^2\io_{18}$ has order $3$ by
Theorem \ref{sph2}(1).
 The two differentials
above this differential
then follow from the extensions. They could also be obtained by the
method of pulling back cycles that we have been using.

When $\nu\le6$, the generator of $E_2^{1,2j}(18,26,34)$ is
$z\equiv h\io_{34}+h^3\io_{26}+h^{\nu+3}\io_{18}$. To obtain the
last term, we used \ref{yang}(3) to write
$h\ot h^2\io_{18}$ as
$d(h^{\nu+3}\io_{18})$ mod lower terms. The leading term of
$\rho\partial_1(z)$ in $E_2^{2,2j}(14)$ is $h^{\nu+3}\ot h\io_{14}$, which is
a generator if $\nu=6$. If $\nu<6$, then this is $d((h^{\nu+2}+L)
\io_{14})$, and so
$z$ pulls back to $z'\equiv z+h^{\nu+2}\io_{14}$. Then $\partial(z')$ in
$E_2^{2,2j}(10)$ is
\begin{equation}\label{4terms}
(h\ot vh^5+h^3\ot vh^3+h^{\nu+3}\ot h^2+h^{\nu+2}\ot h)\io_{10},\end{equation}
using Proposition \ref{coass}.
If $\nu<5$, then the first term is the leading term, and it is a generator.
If $\nu=5$, we must keep track of unit coefficients, since the first and last
terms have the same excess.

Let $j-9=2\cdot3^5c$, with $c\not\equiv0$ mod 3.
We start with $h\io_{34}$. The next term ($h^3\io_{26}$) is insignificant.
The leading term of the image under $\mapright{\partial}E_2^{2,2j}(18)$
is $h\ot\frac12 v^2h^2\equiv \frac12h\ot h^2$. Incorporating coefficients into
the analysis of the previous paragraph, \ref{yang}(3) actually says that
$d(h^8\io_{18})\equiv-3^5ch\ot h^7\io_{18}
\equiv-ch\ot h^2\io_{18}$, and so $z$ is actually
equivalent to $h\io_{34}+h^3\io_{26}+\frac1{2c}h^8\io_{18}$. The leading term
of $\rho\partial_1(z)$ is
$$\tfrac1{2c}h^8\ot(-h)\io_{14}\equiv\tfrac1{2c}h\ot h^6\io_{14}\equiv
-d(\tfrac1{2c}h^7\io_{14}),$$
since by \ref{yang}(3) $d(h^7\io_{14})\equiv-\frac12(2\cdot3^5+2)h\ot h^6
\io_{14}$.
Thus the refined form of $z'$ has significant terms
$h\io_{34}+\frac1{2c}h^7\io_{14}$, and so the leading terms of $\partial(z')$
are $$(h\ot \tfrac14vh^5+\tfrac1{2c}h^7\ot (-h))\io_{10}\equiv
(\tfrac14h\ot h^5+\tfrac1{2c}h\ot h^5)\io_{10},$$
and this is 0 in $E_2$ if $c\equiv1$ mod 3, and is a generator if $c\equiv2$
mod 3.

The boundary into $E_2^{2,2j}(22)$ is $(h\ot h^3+h^3\ot h)\io_{22}$, which is
0 when the group is isomorphic to $\bz/3$.
\end{pf}

The next result also follows by the methods already employed. Note however the
excluded case, which requires major refinements, deferred to the next
section.
\begin{thm}\label{5,8}
If $j$ is odd, and $j\equiv5$ or $8$  mod $9$, but $\nu(j-17)\ne13$, then
$$E_2^{1,2j}(Y_7)\approx E_2^{2,2j}(Y_7)\approx
\bz/3^2\oplus\bz/3^{\min(17,\nu(j-17)+4)}.$$
\end{thm}
\begin{pf} The proof when $j\equiv5$ is particularly simple. The result here
is just that $E_2^{1,2j}(Y_7)\approx
E_2^{2,2j}(Y_7)\approx\bz/3^2\oplus\bz/3^5$.
It is most conveniently seen with charts such as those of
the earlier proofs in this section.
In this case, the two main towers have groups of exponent
2, 1, 1, 1, and $\min(\nu(j-5)+1,5)$, reading from bottom to top. These are the
groups corresponding to generators of dimensions 34, 26, 18, 14, and 10,
respectively. There is
also a group of exponent 2
from the 22-class, and it extends cyclically above the lowest 1.

We will show that
the boundary  is nonzero from enough of the bottom
groups of the $E_2^{1,2j}$-tower to just kill the group of exponent $\min(\nu(j-
5)+1,5)$ at the top of the $E_2^{2,2j}$-tower. That leaves a $\bz/3^5$ in each
tower, and the $\bz/3^2$ coming from the 22-class
cannot be involved in differentials
and must split off for algebraic reasons.

Let $\nu=\nu(j-5)$.
To see these boundary morphisms, we show that the element at the top of
the $E_2^{2,2j}$-tower (i.e., the element of order $3$ in $E_2^{2,2j}(10)$)
is hit by the $\bz/3$ on 14 if $\nu\ge4$,
by the $\bz/3$ on $18$ if $\nu=3$, and by the $\bz/3$ on $26$
if $\nu=2$. Other differentials are seen from the cyclic extensions by
reading down the towers.
The differential when $\nu\ge4$ was proved in \cite[2.4]{BDMi}.
The differential when $\nu=3$ is seen by pulling the generator of
$E_2^{1,2j}(18)$ back to $z\equiv h\io_{18}+h^2\io_{14}$ and then using
\ref{coass} to obtain $\partial(z)\equiv h\ot \frac12v^2h^2\io_{10}$.
The $v^2$ can be moved to the left using \ref{std}(1), and by \ref{sph2}(1)
$\frac12 h\ot h^2\io_{10}$
is an element of order 3 in $E_2^{2,2j}(10)\approx\bz/3^4$.
The case $\nu=2$ is similar, with the leading term of
$\partial(h\io_{26}+h^2\io_{18}+h^4\io_{14})$ being $h\ot vh^3\io_{10}$, which
has order 3 in $E_2^{2,2j}(10)\approx\bz/3^3$. This completes the proof when
$j\equiv5$ mod 9.

Now suppose $j\equiv8$ mod 9, and let $\nu=\nu(j-17)$. The picture is similar
to that just described, with groups of exponent $\min(17,\nu+1)$, 1, 1, 1, and
2, from bottom to top, and a group of exponent 2  extending just
above the lowest 1. These groups correspond to generators of dimensions 34, 26,
18, 14, 10, and 22, as in the case $j\equiv 5$ just considered.
The claim is that differentials from the $E_2^{1,2j}$-tower kill all but the
bottom $3^{17}$ elements in the $E_2^{2,2j}$-tower if $\nu\ge14$, and that
they kill the top $\bz/3^2$ if $\nu\le12$. Actually, when $\nu\ge14$, the
initial element hit also involves a summand in the $22$-summand, but these
elements hit are just the appropriate 3-power times the element at the bottom
of the $E_2^{2,2j}$-tower.
The $E_2^{2,2j}(22)\approx\bz/3^2$ is a split summand in
$E_2^{2,2j}(Y_7)$ even though it may be a summand of a class hit by a
boundary. We will illustrate this carefully in the case $\nu=15$ below.

When $\nu\ge16$, the differential from the bottom of the tower
into the class on the 26-class follows from \cite[2.4]{BDMi}. Of course, the
remaining differentials in this case follow from the extensions.

When $\nu=15$, the generator of $E_2^{1,2j}(34)$ has leading term
$h^{16}\io_{34}$ by \ref{sph}(2), and this pulls back to
$z\equiv h^{16}\io_{34}+uh^{13}\io_{26}$, where $u$ is a unit in $\bz_{(3)}$.
Usually we don't bother to list these unit coefficients, and here the value of
$u$ will not be important, but because cancellation issues will come into play,
we feel that the unit should at least be given lip service. The leading term of
$\partial(z)$ in $E_2^{2,2j}(18)\oplus E_2^{2,2j}(22)$ is
\begin{equation}\label{usum}
uh^{13}\ot 2vh\io_{18}+uh^{13}\ot(-h)\io_{22}.
\end{equation}
Using \ref{yang}, each of
these terms is a generator of its summand. On the other hand, as in the case
$\nu=11$ of the proof of \ref{4}, there is a cycle $z'$ in $C(18,22,26)$
which restricts to a generator of
$E_2^{2,2j}(26)$, and has $3z'$ homologous to a unit times (\ref{usum}).
To clarify the splitting, that $E_2^{2,2j}(22)\approx\bz/3^2$ splits as a
direct summand of $E_2^{2,2j}(Y_7)$, we again use the exact sequence
(\ref{R}). The argument following (\ref{R}) applies verbatim, with
$\bz/3^{13}$ and $\bz/3$ replaced by $\bz/3^{17}$ and $\bz/3^2$, respectively.

The case $\nu=14$ is similar, but involves a 2-step extension process.
In the diagram of the type \ref{charts2}, $E_2^{2,2j}(26)$ extends into
$E_2^{2,2j}(22)\approx\bz/3^2$ and into a $\bz/3^2$ built from $E_2^{2,2j}
(18)$ and $E_2^{2,2j}(14)$.\footnote{$E_2^{2,2j}(10)$ is in the image of
$\partial$, and hence does not figure into the extension question being
considered here.}
The boundary hits into an element of order
3 in each of these summands, which in the case of the second summand means that
it hits a generator of $E_2^{2,2j}(14)$. In order to know that the
splitting is as claimed, we must verify that the element hit is $3^2$ times
a generator of $E_2^{2,2j}(26)$. This is the same sort of verification
that we have been making in some other cases, i.e. that the boundary and the
extension involve classes that are unit multiples of one another, but here the
extension is a 2-step process.

{\bf Boundary:} The generator of $E_2^{1,2j}(34)$ pulls back to
$z\equiv h^{15}\io_{34}+uh^{12}\io_{26}$, with $u$ a unit.
The component of the boundary of this in
$E_2^{2,2j}(22)$ is $uh^{12}\ot(-h)\io_{22}$. On the other hand,
the boundary into $E_2^{2,2j}(18)$ satisfies
$$\partial
(z)\equiv uh^{12}\ot vh\io_{18}\equiv 2uh^{10}\ot h\io_{18}\equiv -2uh\ot
h^8\io_{18}=d(2uh^9\io_{18}).$$
Here we have used the three parts of \ref{yang}, with the last step using that
$2j-18=2(2\cdot 3^{14}c+8)$, and so $\frac14(2j-18)\equiv 1$ mod 3. Thus $z$
pulls back to $z'\equiv z -2uh^9\io_{18}$, and the leading term of
$\partial(z')$ is $2uh^9\ot h\io_{14}$.

{\bf Extension:} Similarly to (\ref{form2}), $d(h^{14})\io_{26}$ is
an element of order 3 in $E_2^{2,2j}(26)$, and it extends to a cycle
$z'$ in $C(18,22,26)$ such that, mod lower classes, $3z'$ is homologous to
$$h^{13}\ot(-h)\io_{22}+h^{13}\ot 2vh\io_{18}.$$
To evaluate $3^2z'$, we use the second 3 to reduce each $h^{13}$ to $h^{12}$.
The second term becomes $$2h^{10}\ot h\io_{18}\equiv -2h\ot h^8\io_{18}
\equiv d(2h^9)\io_{18}\equiv -2h^9\ot(-h)\io_{14}.$$
Here we have applied (\ref{bdryfor}) at the last step.

Thus we have a unit times
$h^{12}\ot(-h)\io_{22}+2h^9\ot h\io_{14}$ as the leading term of
both the image of the boundary, and the $3^2$-multiple of the generator.

The case $\nu\le12$ is much easier. The generator of $E_2^{1,2j}(34)$
pulls back to $z\equiv h^{\nu+1}\io_{34}+h^{\nu-2}\io_{26}+h^{\nu-5}\io_{18}
+h^{\nu-6}\io_{14}$ and this satisfies $\partial(z)\equiv h^{\nu+1}\ot
vh^5\io_{10}$, which is a generator since it does not desuspend.
\end{pf}

The final case differs from the others in that $E_2^{1,2j}(Y_7)$ and
$E_2^{2,2j}(Y_7)$ are not isomorphic.
\begin{thm}\label{2} Assume $j$ is odd and $j\equiv2$ mod $9$.
Then $E_2^{2,2j}(Y_7)\approx\bz/3^2\oplus\bz/3^{\min(13,\nu(j-11)+4)}$, while
$$E_2^{1,2j}(Y_7)\approx\begin{cases}\bz/3^3\oplus\bz/3^5&\text{if $\nu(j-
11)=2$}\\
\bz/3^4\oplus\bz/3^{\min(11,\nu(j-11)+2)}&\text{if $\nu(j-11)>2$.}
\end{cases}$$
\end{thm}
\begin{pf} Let $j$ be as in the theorem, and $\nu=\nu(j-11)$. The picture when
$\nu<8$ is as in Diagram \ref{chart3}.

\begin{diag}\label{chart3}

\begin{center}\begin{picture}(440,190)
\put(117,40){$34$}
\put(117,70){$26$}
\put(132,100){$18$}
\put(102,100){$22$}
\put(132,130){$14$}
\put(132,160){$10$}
\put(125,55){\line(0,1){10}}
\put(140,115){\line(0,1){10}}
\put(140,145){\line(0,1){10}}
\put(120,80){\line(-1,2){10}}
\put(130,80){\line(1,2){10}}
\put(227,42){$2$}
\put(193,105){\frame{$\nu+1$}}
\put(238,105){$\bullet$}
\put(237,164){$2$}
\put(228,72){$\bullet$}
\put(237,132){$\bullet$}
\put(230,48){\line(0,1){22}}
\put(240,108){\line(0,1){22}}
\put(240,139){\line(0,1){24}}
\put(228,80){\line(-1,3){8}}
\put(232,80){\line(1,3){8}}
\put(307,42){$2$}
\put(273,105){\frame{$\nu+1$}}
\put(318,105){$\bullet$}
\put(317,164){$2$}
\put(308,72){$\bullet$}
\put(317,132){$\bullet$}
\put(310,48){\line(0,1){22}}
\put(320,108){\line(0,1){22}}
\put(320,139){\line(0,1){24}}
\put(308,80){\line(-1,3){8}}
\put(312,80){\line(1,3){8}}
\put(306,42){\vector(-1,2){60}}
\put(306,46){\vector(-1,2){60}}
\put(223,25){$E_2^{2,2j}$}
\put(303,25){$E_2^{1,2j}$}
\put(261,5){$\nu<8$}
\end{picture}
\end{center}
\end{diag}

The indicated boundary is seen by pulling back the generator $\a_{m/2}\io_{34}$
to a cycle $z$ on $C(14,18,22,26,34)$, and then obtaining $\a_{m/2}\ot
vh^5\io_{10}$ as the leading term of $\partial(z)$. This generates
$E_2^{2,2j}(10)$. This generator $\a_{m/2}$ is as described in \ref{sph}(1).
The slash does not mean division; this notation was introduced in papers
preceding \cite{BCM}, where it was first applied unstably.

The boundary from $C(26,34)$ into the large group
$E_2^{2,2j}(22)$ has leading term $\a_{m/2}\ot h^3\io_{22}$, which is 0
if $\nu<8$. Here we have $m=\frac12(j-17)$, and we use the argument of
\ref{coass} to see the
factor on the RHS of the $\ot$. (Because of $\a_2$ and $\a_1$ attaching maps,
going from 34 to 22 is like going from
26 to 14, with coefficient $T_4\equiv h^3$ in Proposition \ref{coass}.)
The claimed splitting when $\nu<8$ follows for
algebraic reasons from Diagram \ref{chart3}.

If $\nu=8$, then $\a_{m/2}\ot h^3\io_{22}$ has order 3 in $E_2^{2,2j}
(22)$, by \ref{sph2}(1). If we let $g$ denote a generator of $E_2^{2,2j}
(34)$, then similarly to Diagram \ref{chart3}, $3^3g=a+b$, where $a$ is
detected on the 22-class, and $b$ on the 18-class. We have relations
in $E_2^{2,2j}(Y_7)$ $3^9a$,
$3^4b$, and (from the boundary) $3^2b+3^8ua$, with $u$ a unit in $\bz_{(3)}$.
The quotient group is easily seen to be $\bz/3^{12}\oplus\bz/3^2$, with
generators $g$ and $(1-3^6u)a-3^3g$. The case $\nu=9$ is extremely similar.

If $\nu\ge10$, then $E_2^{2,2j}(22)
\approx\bz/3^{11}$, and the component of the
boundary into this part hits $3^8$ times the generator, as before. But this
implies now that the class on $E_2^{1,2j}(26)$ hits the element of order
3 in $E_2^{2,2j}(22)$. Whereas in the cases $\nu=8$ and 9, the hitting
into the 22-part was without much consequence, because it just
adjoined
another summand to the classes on the 10-cell which were being hit, the boundary
described in the preceding sentence causes one less element in the kernel and
cokernel. In the sort of description given in the previous paragraph, the
relation $3^9a$ is changed to $3^{11}a$. Now we have
$$3^{13}g=3^{10}a=-\tfrac1u3^4b=0,$$
and the claimed splitting follows.
\end{pf}

\section{The final case}\label{lastcase}
In this section, we establish the final and most difficult case of
$E_2^{1,2j}(Y_7)$, with $\nu(j-17)=13$. We will explain why we cannot say for
exactly which such values of $j$ the maximal order is achieved.
\begin{thm}\label{lastthm} If $j$ is odd, and $\nu(j-17)=13$, then for
$\delta$ equal to one of the numbers $2$, $5$, or $8$,
$$E_2^{1,2j}(Y_7)\approx E_2^{2,2j}(Y_7)\approx\bz/3^2\oplus\bz/3^{\min(19,
\nu(j-17-2\delta\cdot3^{13})+4)}.$$
\end{thm}
The methods of this paper do not allow us to determine which of the three
numbers equals $\delta$.
\begin{pf} Let $j-17=2\cdot3^{13}c$, with $c\not\equiv0$ mod 3. The proof
begins just like that of the case $j\equiv8$ mod 9 in Theorem \ref{5,8}.
In the diagram of the type that we have been using, the main tower has groups
of exponent 14, 1, 1, 1, and 2, reading from bottom to top, and a group of
exponent 2 extending above the lowest 1. We choose as
the generator of $E_2^{1,2j}(34)$ the element $-\a_{m/14}\io_{34}$, where
$m=\frac12(j-17)$. We use (\ref{e}) to write it as
$(h^{14}+L)\io_{34}$, with $L$ defined on $S^{27}$. (We choose the minus on
$\a$ to remove the minus signs in (\ref{e}) and (\ref{s}).)

The boundary $E_2^{1,2j}(34)\mapright{\partial}E_2^{2,2j}(26)$
sends this generator to a class congruent mod lower terms to
$$h^{14}\ot 2vh\io_{26}\equiv2h^{12}\ot h\io_{26}\equiv-2h\ot h^{10}\io_{26}
=d((h^{11}+L')\io_{26}).$$
Here we use all three parts of \ref{yang}, with the last step using that
$j-13=2c3^{13}+4$, and so $\frac12(j-13)\equiv 2$ mod a high power of 3.
Thus the generator pulls back to $z\equiv h^{14}\io_{34}-h^{11}\io_{26}$.

Next we consider $\partial(z)$ in both $E_2^{2,2j}(22)$ and in $E_2^{2,2j}
(18)$. The former has leading term $(h^{14}\ot h^3-h^{11}\ot(-
h))\io_{22}$. This desuspends to $S^{19}$ and hence is 0 in $E_2^{2,2j}(
22)\approx\bz/3^2$. Since the 22-cell factor is split from $C(10,14)$,
 we do not need to write this as a boundary and append to $z$.
In $C^2(18)$, we have
$$\partial(z)\equiv (h^{14}\ot\tfrac12v^2h^2-h^{11}\ot 2vh)\io_{18}
\equiv -2h^9\ot h\io_{18}\equiv 2h\ot h^7\io_{18}=d(-\tfrac12h^8\io_{18}),$$
similarly to the previous paragraph. Thus $z$ pulls back to
$z'\equiv z+\frac12h^8\io_{18}$.

The leading term of $\partial(z')$ in $E_2^{2,2j}(14)$ is
$\tfrac12h^8\ot(-h)\io_{14}\equiv \tfrac12h\ot h^6\io_{14}$, and so
$\partial(z')=d((-\frac1{10}h^7+L)\io_{14})$,
since $\frac12(j-7)\equiv5$ mod a high power of 3.
As we will be working at most mod
9, we replace the 10 by 1. Thus $z'$ pulls back to
\begin{equation}\label{z''}
z''\equiv-\a_{m/14}\io_{34}-h^{11}\io_{26}+\tfrac12h^8\io_{18}+h^7\io_{14}.
\end{equation}
Now we use (\ref{s}) for $\a_{m/14}$, and obtain terms in $\partial(z'')$ due
to the first and last terms of (\ref{z''}):
\begin{equation}\label{ch}
\partial(z'')\equiv (ch\ot\tfrac14 vh^5+h^7\ot(-h))\io_{10}\equiv
(\tfrac14c+1)h\ot h^5\io_{10}.\end{equation}
This is a generator if $c\equiv1$ mod 3, in which case the diagram described at
the beginning of the proof has differential from the generator
of $E_2^{2,2j}(34)$ and 3 times the
generator killing $E_2^{2,2j}(10)$, yielding $\bz/3^2\oplus\bz/3^{17}$
as the groups $E_2^{1,2j}(Y_7)$ and $E_2^{2,2j}(Y_7)$,
as claimed in this case. The splitting is true for algebraic reasons.

If $c\equiv2$ mod 3, then $\partial(z'')$ is not a generator of $E_2^{2,2j}
(10)\approx\bz/9$, but it might be 3 times the generator.  This requires
second-order information throughout the entire analysis above. This is
something that we have not had to do in past applications. In particular,
we need finer information in all three parts of Lemma \ref{yang}, in
both descriptions of $\a_{m/e}$ in Theorem \ref{sph}(2), and in Proposition
\ref{coass}.

We now write $c=3k+2$.
The cycle $z''$ above can be written as
\begin{equation}\label{zpp}z''=-\a_{m/14}\io_{34}+(-
h^{11}+A_{10}+L_{10})\io_{26}+(\tfrac12h^8+A_7+L_7)\io_{18}+(h^7+A_6+L_6)
\io_{14},\end{equation}
where $A_i$ has excess exactly $i$, and $L_i$ has excess less than $i$. When we
evaluate $\partial(z'')$, the terms of excess 5 will cancel out as in (\ref{ch})
with $c\equiv2$, and
so we can desuspend $\partial(z'')$ to $S^9$. Our differential into
$E_2^{2,2j}(10)$ is equal to 3 times the generator if and only if
the desuspension of $\partial(z'')$ yields a generator of $E_2^{2,2j-1}(S^9)$.

Let $B_4$ be the terms of excess exactly 4 in $T_3$ of Proposition \ref{coass}.
The terms of
excess 4 or 5 in $\partial(z'')$ are
\begin{equation}\label{parz''}
-\a_{m/14}\ot(\tfrac14vh^5+B_4)+\tfrac12h^8\ot\tfrac12h^2
+(h^7+A_6)\ot(-h).
\end{equation}
Note how certain terms such as $\partial A_{10}$ and $\partial A_7$ were
dropped because they yield terms whose excess is less than 4.

By an analysis similar to \cite[2.11(5)]{BDMi}
we have, when $p=3$ and $c\not\equiv0$ mod 3,
$$\a_{c3^{e-1}/e}\equiv -chv^{c3^{e-1}-1}+\tfrac32ch^2v^{c3^{e-1}-2}
-3ch^3v^{c3^{e-1}-3} \text{ mod }9.$$
Let $h^7\ot(-h)=h\ot h^5+C_4+L_4$, where $C_4$ has excess 4, and $L_4$ excess
less than 4. Omitting terms of excess less than 4, (\ref{parz''}) becomes
\begin{equation}\label{par2}
((3k+2)(h+3h^2)
-3h^2+6h^3)\ot(\tfrac14vh^5+B_4)+\tfrac14h^8\ot h^2+h\ot h^5+C_4
-A_6\ot h,\end{equation}
where the $(h+3h^2)$ comes from $hv^{c3^{e-1}-1}=(v-3h)^{c3^{e-1}-1}h$.
Now write $h\ot vh^5$ as $h\ot h^5-3h^2\ot h^5$.
Using coefficients of 3 to reduce the excess
of terms on the right side of the $\ot$, we can rewrite (\ref{par2}) in excess
4 as
\begin{equation}\label{par3}
\tfrac34kh\ot h^5+\tfrac32(h+\tfrac12h^2+h^3)\ot h^5
+2h\ot B_4+\tfrac14h^8\ot h^2+C_4-A_6\ot h,
\end{equation}
where the $\frac32h\ot h^5$ comes from the $2h\ot\frac14h^5$ and $h\ot h^5$
in (\ref{par2}).
Let $D_4=\tfrac32(h+\tfrac12h^2+h^3)\ot h^5+\tfrac14h^8\ot h^2+C_4$, a
specific class of excess 4, independent of the value of $k$ and of any choices
of the sort that we are about to mention.

The term $B_4$ is the terms of excess 4 in (\ref{T32}). It could also have
included
any terms of excess 4 in the homogeneous part of $T_3$ discussed in the
paragraph after (\ref{T32}), but as discussed there, this homogeneous part
has excess
less than 4. Then $B_4$ contains a term $-\frac12v^2h^4$ which appears in
(\ref{T32}), and it could contain a term $9c_2h^6$ if $c_2\not\equiv0$ mod 3.
However, because of a term with coefficient $\tfrac{59}6$ which has $c_2$ as
coefficient, we can infer that $c_2\equiv0$ mod 3. Thus $B_4=-\frac12v^2h^4$,
and so we can let $D_4'=D_4-h\ot v^2h^4$, still a specific element of excess 4,
and we have
\begin{equation}\label{par4}
\tfrac14kh\ot h^4+D_4'-A_6\ot h
\end{equation}
as our new expression for $\partial(z'')$ mod $L$.

Next we study $A_6$. To find it, we apply $\mapright{\partial}E_2^{2,2j}(14)$
to the sum $z_3=X_1+X_2+X_3$ of the first three terms
of (\ref{zpp}), and write the result as $d(A_6)$.
The terms in $\partial(X_1)$ will have excess less
than 5, and so may be omitted from the analysis. There is one term,
$\frac12h^8\ot(-h)$, of excess 6, which accounts for the $h^7$ in (\ref{z''}).
There are a number of terms of excess 5, which contribute toward $A_6$.
In particular, note that $d(h^6)\equiv
h\ot h^5$, and so each occurrence of $h\ot
h^5$ in $\partial(z_3)$ affects the coefficient of $h^6$ in $A_6$.
The leading part of $\partial(X_2)$ is $-h^{11}\ot T_4$, where $T_4$ is
as in \ref{coass}. The full form of $T_4$ is given in (\ref{T4}) and involves
a homogeneous part whose coefficient $c$ we do not know. Two parts of this
homogeneous part have a factor of 3, which can be used to reduce the excess,
but $ch^{11}\ot v^2h\equiv ch^7\ot h\equiv -ch\ot h^5$ will cause a
$ch^6$-term in $A_6$, and hence a $ch\ot h^4$ in (\ref{par4}).
Thus the coefficient of $h\ot h^4$ in (\ref{par4}) is $k+D+c\in \bz/3$, where
$D$ is something which we could compute if we really needed to. Note also that
for our purposes (\ref{par4}) lies in $\bz/3$ generated by $h\ot h^4$.
The coefficient $c$ has a value; we just don't know how to find it.
Therefore, there is one value of $k$ in $\bz/3$ for which (\ref{par4}) is 0.
(The diligent reader can check that such considerations cannot affect
earlier parts of the argument.) Thus the differential into $E_2^{2,2j}(10)$
 is 0 if and only if $k$, defined by $j-17=2(3k+2)3^{13}$, has this
value mod 3. Letting $\delta=3k+2$ mod 9, this establishes the theorem.
\end{pf}

\section{Periodic homotopy of $E_7$}\label{E7}
In this section we use the results for $E_2^{s,2j}(Y_7)$
already achieved to deduce
that $v_*(E_7)$ is as claimed in Theorem \ref{main}. The first result
almost finalizes $v_*(Y_7)$, given the results for $E_2^{s,2j}(Y_7)$ determined
in the previous two sections.
\begin{thm}\label{x7} The $v_1$-periodic UNSS of $Y_7$ converges to $v_*(Y_7)$.
If $j$ is odd, then $v_{2j+1}(Y_7)=0$,
$v_{2j}(Y_7)\approx
v_{2j}(S^7)$, $v_{2j-2}(Y_7)\approx E_2^{2,2j}(Y_7)$, and there is
an exact sequence
$$0\to v_{2j-1}(S^7)\to v_{2j-1}(Y_7)\to E_2^{1,2j}(Y_7)\to 0.$$
\end{thm}
\begin{pf} The main thing that we have to worry about in proving convergence of
the $v_1$-periodic UNSS is to rule out the possibility of a $v_1$-periodic
homotopy class which is not seen in $v_1$-periodic $E_2$. This could come about
by having a sequence of homotopy classes related by a filtration-increasing
$v_1$-multiplication. Such a class could also be the target of a \lq\lq
differential'' from an element of $v_1$-periodic $E_2$. The way that we will
show that these things cannot happen for $Y_7$ is to note that $Y_7$ is built
by fibrations from spaces where we have already established convergence.

In (\ref{split12}), it was noted how the
$v_1$-periodic UNSS of $Y_7$ splits into
the part from $S^7$ and the part from even-dimensional classes. As all of this
is confined to filtrations 1 and 2, we obtain the following schematic picture
for $E_2^{s,t}(Y_7)$, which must necessarily equal $E_\infty$.

\begin{center}
\begin{tabular}{c|c|c|c|c|c}
\cline{2-5}
$s=2$&\ ev\ &\ $S^7$\ &&\qquad&\\
\cline{2-5}
$s=1$&&\ ev\ &\ $S^7$\ &&\\
\cline{2-5}
$t-s=$&$2j-2$&$2j-1$&$2j$&$2j+1$&$j$ odd
\end{tabular}
\end{center}

Here a box labeled $S^7$ means the
corresponding group $E_2^{s,t}(S^7)$, while a box
labeled \lq\lq ev'' (for \lq\lq even'') means the corresponding group
$E_2^{s,t}(10,14,18,22,26,34)$, as computed in Section \ref{Y7}. This $E_2$
calculation is consistent with the fibrations (\ref{X7fib}) and
$S^7\to\Om W\to\Om S^{23}$ of Proposition \ref{HarK}.

For $X=\Om S^{23}$, $\Om B(11,15)$, or $\Om E_7/F_4$, the $v_1$-periodic UNSS
collapses to isomorphisms, if $j$ is odd,
$$\vp_{2j+\eps}(X)\approx\begin{cases}0&\text{if $\eps=0$ or 1}\\
E_2^{2,2j}(X)&\text{if $\eps=-2$}\\
E_2^{1,2j}(X)&\text{if $\eps=-1$.}\end{cases}$$
This is true for $\Om S^{23}$ by \cite[6.1]{BCR}, for $\Om B(11,15)$ by the
fibration $\Om S^{11}\to \Om B(11,15)\to\Om S^{15}$, and for $\Om E_7/F_4$
by Theorem \ref{BTthm}. (Although \ref{BTthm} dealt with convergence for
$E_7/F_4$, the methods of Section \ref{Y7} show that the calculation for
$E_2(\Om E_7/F_4)$ is just that for $E_2(E_7/F_4)$ shifted back by 1 dimension,
and of course the same is true of $v_1$-periodic homotopy groups.)

Let $j$ be odd.
We can use a Five Lemma argument once we establish that, for
$\eps=1$ or 2, there are morphisms
$v_{2j-\eps}(-)\to E_2^{\eps,2j}(-)$ for these spaces. To see that such
morphisms exist, we note that since compact Lie groups and spheres have
$H$-space exponents (\cite{James}),
the spaces with which we deal here all have $H$-space exponents. By \cite{DM},
this implies that each $v_1$-periodic homotopy group is a direct summand of
some actual homotopy group, and then we can take the morphism from homotopy
to homotopy mod filtration greater than $\eps$, which is (unlocalized)
$E_\infty^\eps$, then to (unlocalized)
$E_2^\eps$ as the kernel of the differentials, and then to $v_1$-periodic
$E_2^\eps$. This argument is similar to that used in \cite{D}.

Thus, letting $X=E_7/F_4$ and $B=B(11,15)$, there is a commutative diagram of
exact sequences
$$\begin{CD}
v_{2j-1}(\Om X)@>>>v_{2j-2}(\Om S^{23}\times \Om B)@>>>v_{2j-2}(Y_7)@>>>v_{2j-
2}(\Om X)@>>>0\\
@V\approx VV @V\approx VV @VVV@V\approx VV\\
E_2^{1,2j}(\Om X)@>>>E_2^{2,2j}(\Om S^{23}\times\Om B)@>>>E_2^{2,2j}(Y_7)
@>>>E_2^{2,2j}(\Om X)@>>>0
\end{CD}$$
which implies that $v_{2j-2}(Y_7)\to E_2^{2,2j}(Y_7)$ is an isomorphism.

Similarly, there is a commutative diagram with exact rows and the first column
exact
$$\begin{CD}@.0@.@.@.\\
@.@VVV@.@.@.\\
@.E_2^{2,2j+1}(S^7)@>\approx>>E_2^{2,2j+1}(Y_7)@.@.\\
@.@VVV@VVV@.@.\\
0@>>>v_{2j-1}(\Om W\times\Om B)@>>>v_{2j-1}(Y_7)@>>>v_{2j-1}(\Om X)@>>>
v_{2j-2}(\Om W\times \Om B)\\
@.@VVV@VVV@V\approx VV@V\approx VV\\
0@>>>E_2^{1,2j}(\Om S^{23}\times\Om B)@>>>E_2^{1,2j}(Y_7)@>>>E_2^{1,2j}(\Om X)
@>>>E_2^{2,2j}(\Om S^{23}\times\Om B)\\
@.@VVV@.@.@.\\
@.0@.@.@.   \end{CD}$$
which implies that the second column fits into a short exact sequence.

The portion of the theorem about $v_{2j+1}(Y_7)$ and $v_{2j}(Y_7)$ is immediate
from the exact sequence in $v_*(-)$ associated to the fibration (\ref{X7fib}).
\end{pf}

We restate the following result from \cite[1.3(1)]{BDMi}.
\begin{lem}\label{BS} The projection map $B(3,7)\to S^7$ induces an isomorphism
in $v_{2j-1}(-)$ unless $j$ is odd and $j\equiv21$ mod $27$,
in which case it is
a surjection $\bz/3^4\to\bz/3^3$. The isomorphic groups are $0$ if $j$ is even,
while if $j$ is odd, they are cyclic of order
$3^{\min(3,1+\nu(j-3))}$.
\end{lem}

The next result, combined with the above results and Theorems \ref{1,7},
\ref{4}, \ref{0}, \ref{5,8}, and \ref{lastthm}
gives $v_*(E_7)$ for most values of $*$.
\begin{thm} $(a.)$ If
$j$ is odd, $j\not\equiv2$ mod $9$, and $j\not\equiv21$ mod $27$,
then the exact sequence of the fibration $Y_7\to B(3,7)\to E_7$ breaks up into
isomorphisms $$v_{2j}(Y_7)\mapright{\approx}v_{2j}(B(3,7))\text{\quad and\quad}
v_{2j-1}(E_7)\mapright{\approx}v_{2j-2}(Y_7)$$ and a short exact sequence
$$0\to v_{2j}(E_7)\to v_{2j-1}(Y_7)\mapright{\phi}v_{2j-1}(B(3,7))\to 0.$$
If $E_2^{1,2j}(Y_7)\approx\bz/3^{e_1}\oplus\bz/3^m$, with $1\le e_1\le2$, is as
given in Theorems \ref{1,7}, \ref{4}, \ref{0}, \ref{5,8}, and \ref{lastthm},
and
$v_{2j-1}S^7\approx v_{2j-1}(B(3,7))\approx\bz/3^{e_2}$ is as in \ref{BS}, then
\begin{equation}\label{spl}
v_{2j-1}(Y_7)\approx\bz/3^{e_1+e_2}\oplus\bz/3^m,\end{equation}
and $\phi$ sends the first summand onto $\bz/3^{e_2}$.

$(b.)$ If $j$ is even, then $v_{2j}E_7= v_{2j-1}E_7=0$.\label{not2}
\end{thm}
Note that even if $\phi$ sent the second summand nontrivially, its kernel
would still be $\bz/3^{e_1}\oplus\bz/3^m$, since $m\ge e_1+e_2$. Thus if
$j$ is as in Theorem \ref{not2}(a.), there are abstract isomorphisms
$v_{2j}(E_7)\approx E_2^{1,2j}(Y_7)$ and $v_{2j-1}(E_7)\approx E_2^{2,2j}(Y_7)$,
with $E_2^{s,2j}(Y_7)$ as given in Theorems \ref{1,7}, \ref{4}, \ref{0},
\ref{5,8}, and \ref{lastthm}. This implies Theorem \ref{main} in these cases.
\begin{pf}
There is a commutative diagram of fibrations
\begin{eqnarray}
\Om W&\to B(3,7)\to&K\nonumber\\
\da\quad&\mapdown{=}&\ \da\label{l0}\\
Y_7\quad &\to B(3,7)\to&E_7\nonumber
\end{eqnarray}
where the last map is the composite $K\to F_4\to E_7$.
Since by \cite[2.10(i)]{F4} the composite $S^7\to \Om W\to B(3,7)\to S^7$ has
degree $3$, we deduce the same of the composite $S^7\to Y_7\to B(3,7)\to S^7$.
We already know that $v_{2j}S^7\to v_{2j}Y_7$ is an isomorphism, and
$v_{2j}B(3,7)\to v_{2j}S^7$ is multiplication by $3$ on isomorphic groups.
It follows that $v_{2j}Y_7\to v_{2j}B(3,7)$ is an isomorphism.

There is a commutative diagram of fibrations
\begin{eqnarray}
 S^7&\to \Om W\to&\Om S^{23}\nonumber\\
\da\ &\da&\label{l1}\\
 S^7&\to Y_7\quad\ &\nonumber.
\end{eqnarray}
The cyclic extension in $v_{2j-1}(\Om W)$ was established in
\cite[pp.294-5]{F4}. This implies the nontrivial extension in $v_{2j-1}Y_7$
claimed in the theorem from the $\bz/3^{e_1}$ on the
22-class in $E_2^{1,2j}(Y_7)$ to $v_{2j-1}(S^7)$ in the exact sequence of
Lemma \ref{x7}.

There cannot be an
extension in $v_{2j-1}Y_7$ from the $\bz/3^m$-summand of $E_2^{1,2j}Y_7$
because of the splitting $F_4=K\times B(11,15)$. The element of order $3$ in
the large summand of $v_{2j-1}Y_7$ comes from $B(11,15)$, while the $S^7$ lies
in $K$. This is made explicit in the commutative diagram of fibrations
\begin{eqnarray}
\Om B(11,15)\times\Om W&\to B(3,7)\to&F_4\nonumber\\
\da\qquad&\da&\ \da\label{l2}\\
Y_7\qquad&\to B(3,7)\to&E_7\nonumber
\end{eqnarray}
That $\phi$ sends the first summand of (\ref{spl})
onto $v_{2j-1}B(3,7)$ follows from the
diagram (\ref{l2}) and the surjectivity of $v_{2j-1}(\Om W)\to v_{2j-1}B(3,7)$
established in \cite[pp.297-8]{F4}.
\end{pf}

One of the cases omitted in the previous theorem is covered in the following
result, the proof of which is very similar.
\begin{thm} If $j$ is odd and $j\equiv21$ mod $27$, then the exact sequence
$($with $B=B(3,7))$
$$0\to v_{2j}Y_7\mapright{\phi_1}v_{2j}B\to v_{2j}E_7\to v_{2j-1}Y_7
\mapright{\phi_2}v_{2j-1}B\to v_{2j-1}E_7\to v_{2j-2}Y_7\to 0$$
has $\phi_1$ an injection $\bz/3^3\hookrightarrow\bz/3^4$, and $\phi_2$ a
surjection from the first summand in $\bz/3^4\oplus\bz/3^m\to\bz/3^4$.
Moreover, $$v_{2j}E_7\approx\coker\phi_1\oplus\ker\phi_2\approx
\bz/3\oplus\bz/3^m.$$
\end{thm}
\begin{pf} Similarly to the previous proof,
the morphism $\phi_1$ follows from \cite[2.5]{BDMi} and \cite[2.10(i)]{F4},
the structure of $v_{2j-1}Y_7$ follows from
(\ref{l1}), and the morphism $\phi_2$ follows
from (\ref{l2}). The $\bz/3^m$ in $v_{2j-1}Y_7$
cannot extend cyclically with $\coker\phi_1$ in $v_{2j}E_7$ because the element
of order $3$ in $\bz/3^m$ lies in $v_{2j-1}\Om B(11,15)$, while $\coker\phi_1$
lies in $v_{2j}K$, and these cannot be related by a $\cdot 3$-extension
due to the splitting $F_4=K\times B(11,15)$.
\end{pf}

We begin working toward determination of $v_{2j-\eps}E_7$ when $j\equiv2$ mod 9
with the following proposition.
\begin{prop} \label{3items}If $j$ is odd and $j\equiv2$ mod $9$,
then the exact sequence of the
fibration $Y_7\to B(3,7)\to E_7$ yields
\begin{itemize}
\item $v_{2j}Y_7\to v_{2j}B(3,7)$ is an isomorphism of $\bz/3$'s;
\item $v_{2j-1}E_7\to v_{2j-2}Y_7$ is an isomorphism;
\item $v_{2j}E_7\approx\ker(v_{2j-1}Y_7\mapright{\phi}v_{2j-
1}B(3,7)\approx\bz/3)$.
\end{itemize}
\end{prop}
\begin{pf} Surjectivity of $v_{2j-1}Y_7\to v_{2j-1}B(3,7)$ follows from
(\ref{l2}), while $v_{2j}Y_7\to v_{2j}B(3,7)$ is bijective as in the proof of
\ref{not2}.\end{pf}
By \ref{3items}, \ref{x7}, and \ref{2}, $v_{2j-1}E_7$ is seen to be as claimed
in Theorem \ref{main} when $j\equiv 2$.
It remains to determine $v_{2j-1}Y_7$ and $\phi$, from which $v_{2j}E_7$
follows.
\begin{thm} Let $j$ be odd, and $\nu=\nu(j-11)$. If $2\le\nu\le9$, then
$$v_{2j-1}Y_7\approx\bz/3^{\nu+3}\oplus\bz/3^4$$
and $\phi$ sends $\bz/3^4$ nontrivially. Thus
$\ker\phi\approx\bz/3^{\nu+3}\oplus\bz/3^3$, regardless of
$\phi\bigm|\bz/3^{\nu+3}$.\label{penult}
\end{thm}
\begin{pf} Similarly to the proof of Theorem \ref{not2}, the extension in
$$\bz/3\approx v_{2j-1}S^7\to v_{2j-1}Y_7\to E_2^{1,2j}Y_7\approx
\bz/3^{\nu+2}\oplus\bz/3^{4}$$
is nontrivial from the first summand. From \cite[2.12]{F4}, $v_{2j-1}(\Om
W)\to v_{2j-1}B(3,7)$ is a surjection $\bz/3^{\nu+2}\to\bz/3$, and from
(\ref{l0}) it factors as
$$v_{2j-1}(\Om W)\to v_{2j-1}Y_7\mapright{\phi}v_{2j-1}B(3,7).$$
From (\ref{l1}), $v_{2j-1}(\Om W)\to v_{2j-1}Y_7$
is an injection $\bz/3^{\nu+2}\to\bz/3^{\nu+3}\oplus\bz/3^4$,
since the element of order 3 in $v_{2j-1}(\Om W)$, which comes from
$v_{2j-1}(S^7)$, maps nontrivially. The result now
follows from elementary algebra.
\end{pf}
The same ingredients imply the following result.
\begin{thm} If $j$ is odd and $\nu(j-11)\ge 10$, but
$j\not\equiv11+2\cdot3^{10}$ mod $2\cdot3^{11}$, then
$$v_{2j-1}Y_7\approx\bz/3^{12}\oplus\bz/3^4$$
and $\phi$ is surjective in Proposition \ref{3items}.
\label{cannot}\end{thm}

We cannot deduce from this which summand(s) of $v_{2j-1}Y_7$ maps nontrivially
under $\phi$, and so we cannot tell whether $\ker\phi$ is
$\bz/3^{12}\oplus\bz/3^3$ or $\bz/3^{11}\oplus\bz/3^4$. We suspect that
$\bz/3^4$ maps across, which would imply the first splitting.

Finally we have the following result in the exceptional case. In order to keep
the statement of Theorem \ref{main} readable, we did not distinguish there
between this case, in which we know the precise structure of $v_{2j-1}E_7$,
and the case of Theorem \ref{cannot}, where we do not.

\begin{thm}If $j\equiv11+2\cdot3^{10}$ mod $2\cdot3^{11}$, then
$$v_{2j}E_7\approx\bz/3^{12}\oplus\bz/3^3.$$
\end{thm}
\begin{pf} As in the proof of \ref{penult}, $v_{2j-
1}Y_7\approx\bz/3^{12}\oplus\bz/3^4$. It was shown in \cite[2.12]{F4}  that
$v_{2j-1}(\Om W)\to v_{2j-1}B(3,7)$ is 0 if $j\equiv11+2\cdot3^{10}$ mod
$2\cdot 3^{11}$.

Let $G$ denote the fiber of $K\to E_7$. There is a commutative diagram of
fibrations
$$\begin{array}{ccc}
\Om W&\mapright{=}\ \Om W\qquad&\\
\da&\da&\\
Y_7&\to B(3,7)\to&E_7\\
\da&\da&\da\\
G&\to\quad K\quad\to &E_7.
\end{array}$$
It follows from the Serre spectral sequence of the fibration $\Om W\to Y_7\to
G$ that
$$BP_*(G)\approx BP_*[x_{10},x_{14},x_{18},x_{26},x_{34}],$$
and so charts for $v_*G$ are like charts for $E_2^{s,2j}(Y_7)$
without the part on the
22-class. The chart for $v_{2j-1}G$ and $v_{2j-2}G$ whenever $j\equiv2$ mod 9
is like Diagram \ref{chart3} without the $\nu+1$. In particular, $v_{2j-1}G$
is cyclic with generator on the 26-class. The proof of Theorem \ref{2} in the
case $\nu\ge10$, where it says that the class on $E_2^{1,2j}(26)$ hits the
element of order $3$
in $E_2^{2,2j}(22)$, implies that $v_{2j-1}G\to v_{2j-1}K$
sends the generator to the element of order $3^2$. Now it follows from the
following commutative diagram with exact rows that $\phi$ is surjective on the
$\bz/3^4$ summand.
$$\begin{array}{cccc}
\bz/3^{11}&\bz/3^{12}\oplus\bz/3^4&\bz/3^5&\\
0\to v_{2j}W\to&v_{2j-1}Y_7&\to v_{2j-1}G\to&v_{2j-1}W\\
&\mapdown{\phi}&\da&\mapdown{=}\\
&v_{2j-1}B&\to v_{2j-1}K\to&v_{2j-1}W\\
&\bz/3&\bz/3^{12}&\bz/3^{11}
\end{array}$$
\end{pf}

\def\line{\rule{.6in}{.6pt}}

\end{document}